\documentclass[11pt]{article}
\usepackage{amssymb,amsthm}
\usepackage{amsmath}
\usepackage{amscd}
\usepackage{epsf,graphicx}
\usepackage[all]{xy}
\newtheorem{thm}{Theorem}
\newtheorem{prop}[thm]{Proposition} 

\newtheorem{defi}[thm]{Definition}
\newtheorem{exa}[thm]{Example}

\newcommand{\tr}{{\rm{tr}}}
\newcommand{\ad}{{\rm{ad}}}
\newcommand{\ob}{\rm{Ob}}
\newcommand{\gr}{{\rm{gr}}}

\newcommand{\Der}{{\rm{Der}}}

\newcommand{\End}{{\rm{End}}}

\begin{document} \setlength{\baselineskip}{16pt}
\title{Graphical introduction to classical Lie algebras}
\author{Rafael D\'\i az\thanks{Work partially supported by UCV.}\ \
and Eddy Pariguan\thanks{Work partially supported by FONACIT.}}
\maketitle
\begin{abstract}

We develop a graphical notation to introduce classical Lie algebras.
Although this paper deals with well-known results, our pictorial
point of view is slightly different to the traditional one. Our
graphical notation is fairly elementary and easy to handle and, thus
provides an effective tool for computations with classical Lie
algebras. More over, it may be regarded as a first and foundational
step in the process of uncovering the categorical meaning of Lie
algebras.

\end{abstract}
\section{Introduction}

A first step in the study of an arbitrary category $\mathcal{C}$ is
to define the appropriated set $\mathcal{S(C)}$ of isomorphisms
classes of {\em simple} objects in $\mathcal{C}$. For example, an
object $y \in \ob(\mathcal{C})$ of an abelian category $\mathcal{C}$
(see \cite{Man}) is said to be simple if in any exact sequence
$$0 \to x \to y \to z \to 0,$$
either $x$ is isomorphic to $0$ or $z$ is isomorphic to $0$. It is a
remarkable fact that non-equivalent categories may very well have
equivalent sets of isomorphisms classes of simple objects. Let us
introduce a list of categories that at first seem to be utterly
unrelated and yet the corresponding sets of simple objects are
deeply connected.

\begin{enumerate}
\item $\mathbf{WeylGroup} \subset\mathbf{FinGroup} \subset
\mathbf{LieGroup} \subset \mathbf{Group}.$

We denote by $\mathbf{Group}$ the category whose objects are
groups and whose morphisms are group homomorphisms. We let
$\mathbf{LieGroup}$ denote the subcategory of $\mathbf{Group}$
whose objects are finite dimensional complex Lie groups. Morphism
in $\mathbf{LieGroup}$ are smooth group homomorphisms. We define
$\mathbf{FinGroup}$ to be the full subcategory of $\mathbf{Group}$
whose objects are finite groups and $\mathbf{WeylGroup}$ to be the
full subcategory of $\mathbf{Group}$ whose objects are Weyl
groups.

$S(\mathbf{Group})$ denotes the set of isomorphisms classes of
groups having no proper normal subgroups. $S(\mathbf{LieGroup})$
denotes the set of isomorphisms classes of  Lie groups which are
simple as groups and also are connected and simply connected.
$S(\mathbf{FinGroup})$ denotes the set of isomorphisms classes of
finite simple groups. $S(\mathbf{WeylGroup})$ denotes the set
isomorphisms classes of Weyl groups, which can be taken to be
$A_n=S_{n+1}$, $B_n=C_n=\mathbb{Z}_2^{n}\ltimes S_n$ ,
$D_n=\mathbb{Z}_2^{n-1}\ltimes S_n$, where $S_n$ is the group of
permutations in $n$ letters, and $E_6, E_7, E_8, F_4, G_2$ are so
called exceptional groups.

\item $\mathbf{LieAlg}$ denotes the category whose objects are
finite dimensional complex Lie algebras, morphism are Lie algebra
homomorphism. $S(\mathbf{LieAlg})$ is the set of isomorphisms
classes of simple Lie algebras, i.e., Lie algebras having no proper
ideals.

\item $\mathbf{Root}$ denotes the category of root systems.

\begin{itemize}
\item Objects in $\mathbf{Root}$ are triples $(V,\langle \mbox{ }, \mbox{ }\rangle, \Phi)$ such that

\begin{itemize}

\item{$(V, \langle \mbox{ }, \mbox{ }\rangle)$ is an Euclidean
space.} \item{$\Phi$ is a finite set spanning $V$.} \item{If
$\alpha\in\Phi$ then $-\alpha\in\Phi$, but $k \alpha\not\in \Phi$ if
$k$ is any real number other then $\pm 1$.} \item{For
$\alpha\in\Phi$ the reflection $S_{\alpha}$ in the hyperplane
$\alpha^{\perp}=\{x\in V: \langle x,\alpha\rangle=0\}$ maps $\Phi$
to itself.} \item{For $\alpha, \beta \in\Phi$, $A_{\alpha, \beta}=
{\displaystyle 2\frac{\langle \alpha, \beta\rangle}{ \langle
\alpha, \alpha \rangle}\in \mathbb{Z}}$.}
\end{itemize}

\item Morphisms in $\mathbf{Root}$ from $(V_1,\langle \mbox{ }, \mbox{ }\rangle_1, \Phi_1)$ to
 $(V_2,\langle \mbox{ }, \mbox{ }\rangle _2, \Phi_2)$ consists of
linear transformations $T:V_1 \longrightarrow V_2$ such that
$\langle T(x),T(y)\rangle _2= \langle x,y \rangle _1$ for all
$x,y\in V_1$, and $T(\Phi_1)\subset \Phi_2$.

\item Suppose that $(V_i,\langle \mbox{ }, \mbox{ }\rangle _{V_i}, \Phi_i)$,
$i=1,\dots,n$ are root systems, then the direct sum\\
$V=\bigoplus_{i=1}^{n} V_i$ is an Euclidean vector space with
inner product
$$\langle \mbox{ },\mbox{ }\rangle_V= \sum_{i=1}^{n}\langle \mbox{ },\mbox{ } \rangle _{V_i}$$
Moreover setting $\Phi=\bigsqcup_{i=1}^{n} \Phi_i$, the triple
$(V,\langle \mbox{ }, \mbox{ }\rangle _V, \Phi)$ becomes a root
system called the direct sum of root systems $(V_i,\langle \mbox{
}, \mbox{ }\rangle_{V_i}, \Phi_i)$. $S(\mathbf{Root})$ is the set
of isomorphisms classes of simple root systems, i.e., root systems
which are not isomorphic to the direct sum of two non-vanishing
root systems.

\end{itemize}

\item $\mathbf{Dynkin}$ denotes the category of Dynkin diagrams.
\begin{itemize}
\item Objects in $\mathbf{Dynkin}$ are called Dynkin diagrams and
are non-directed graphs $\Delta$ with the following properties

\begin{itemize}
\item{The set $V_{\Delta}$ of vertices of $\Delta$ is equal to
$\{1,\dots,n\}$ for some $n\geq 1$.} \item{The number of edges
joining two vertices in $\Delta$ is $0,1,2$ or $3$.}
\item{If vertices $i$ and $j$ are joined by $2$ or $3$ edges, then an arrow is chosen
pointing either from $i$ to $j$, or from $j$ to $i$.}
\item{The quadratic form
$$Q(x_{1},x_{2}, \cdots, x_{n})=2\sum_{i=1}^{n}x_{i}^{2}-
\sum_{i\neq j}\sqrt{n_{ij}}x_{i}x_{j}$$ where $(n_{ij})$ is the
adjacency matrix of $\Delta$, i.e., $n_{ij}$ equal the number of
edges from vertex $i$ to vertex $j$ is positive definite.}
\end{itemize}
\item Morphism in $\mathbf{Dynkin}$ from diagram $\Delta_1$ to
diagram $\Delta_2$ consists of maps $\rho: V_{\Delta_1}\rightarrow
V_{\Delta_2}$ such that $Q_2(x_{\rho(1)},x_{\rho(2)}, \cdots,
x_{\rho(n)})=Q_1(x_{1},x_{2}, \cdots, x_{n}).$
\item{
$S(\mathbf{Dynkin})$ denotes the set of of isomorphisms classes of
connected Dynkin diagrams}

\end{itemize}

\end{enumerate}

\noindent Next theorem gives an explicit characterization of
$S(\mathbf{Dynkin})$.
\begin{thm}$S(\mathbf{Dynkin})$ consists of the Dynkin diagrams
included in the following list
\begin{figure}[!h]
\begin{center}
\includegraphics[width=2in]{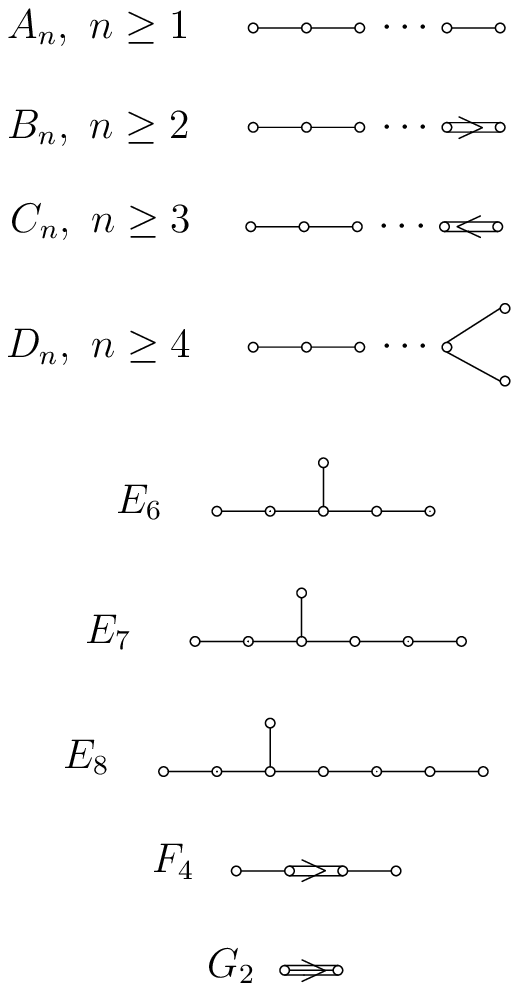}
\caption{Simple Dynkin diagrams. \label{dynkin}}
\end{center}
\end{figure}
\end{thm}
\newpage

\noindent We enunciate the following fundamental

\begin{thm}\label{t1}
\begin{enumerate}
\item{$\mathbf{WeylGroup}\subset S(\mathbf{FinGroup})\subset
S(\mathbf{LieGroup})\subset S(\mathbf{Group})$.}
\item{$S(\mathbf{LieGroup})\cong S(\mathbf{LieAlg}) \cong
S(\mathbf{Root})\cong
S(\mathbf{Dynkin})\rightarrow\mathbf{WeylGroup}$. }
\end{enumerate}
\end{thm}

\noindent Part 1 of Theorem \ref{t1} is obvious. Although we shall not give a
complete proof of part 2 the reader will find in the body of this
paper many statements that shed light into its meaning. The map
$S(\mathbf{Dynkin})\longrightarrow\mathbf{WeylGroup} $ is surjective
but if fails to be injective. Diagrams $B_n$ and $C_n$ of the list
above have both $\mathbb{Z}_{2}^{n}\rtimes S_n$ as associated Weyl
groups.
\section{Lie Algebras}

We proceed to consider in details the category of Lie algebras.
First we recall the notion of a Lie group.

\begin{defi}A group $(G,m)$ is said to be a Lie group if
\begin{enumerate}
\item{$G$ is a finite dimensional smooth manifold.} \item{ the map
$m:G \times G \longrightarrow  G$ given by $m(a,b)=ab$, for all $a,b
\in G$, is smooth.}
\item{ The map $I:G \to G$ given by
$I(a)=a^{-1}$, for all $a \in G$, is smooth.}
\end{enumerate}

\end{defi}

\begin{defi}
A Lie algebra $({\mathfrak g}, [\mbox{ },\mbox{ }])$ over a field
$k$ is a vector space ${\mathfrak g}$ together with a binary
operation $[ \mbox{} \ ,\mbox{}\ ]:{\mathfrak g} \times {\mathfrak
g} \longrightarrow {\mathfrak g}$, called the Lie bracket,
satisfying

\begin{enumerate}
\item{$[ \mbox{} \ ,\mbox{}\ ]$ is a bilinear operation.}
\item{Antisymmetry: $[x,y]=-[y,x]\ \mbox{ for each }\ x,y\in
{\mathfrak g}$.} \item{Jacobi identity: $[x,[y,z]]=[[x,y],z]+
[y,[x,z]]\ \mbox{ for each }\ x,y,z\in {\mathfrak g}$.}
\end{enumerate}

\end{defi}

\begin{exa}
A $k$-algebra $A$ may be regarded as a Lie algebra $(A,[\mbox{
},\mbox{ }])$, with bracket $[x,y]=xy-yx$ for all $x,y \in A$. In
particular $\End(V)$ is a Lie algebra for any $k$-vector space
$V$.
\end{exa}

\begin{exa}
Let $M$ be a smooth manifold. The space
$\Gamma(M)=\{X:M\longrightarrow TM,\ X(m)\in T_m M,\  m\in M\}$ of
vector fields on $M$ is a Lie algebra with the Lie bracket
$$[X,Y]=X^{j}\frac{\partial
Y^{i}}{\partial x^{j}}\frac{\partial}{\partial x^{i}}-
Y^{j}\frac{\partial X^{i}}{\partial x^{j}}\frac{\partial}{\partial
x^{i}},  \ \ \mbox{for all}\ \ X,Y\in\Gamma(M).$$
\end{exa}

\begin{exa}
Let $G$ be a Lie group. The space $T_{e}(G)$ tangent to the identity
$e\in G$ is a Lie algebra since $T_{e}(G)\cong \Gamma(G)^G$ is a Lie
subalgebra of $\Gamma(G)$.

\end{exa}

\begin{defi}
A {\em morphism} of Lie algebras $\rho: {\mathfrak g}
\longrightarrow {\mathfrak h}$ is a linear map $\rho$ from
${\mathfrak g}$ to ${\mathfrak h}$ such that
$\rho([x,y])=[\rho(x),\rho(y)]$ for $x,y \in {\mathfrak g}.$ A
{\em representation} $\rho$ of a Lie algebra $\mathfrak g$ on a
$k$-vector space $V$ is a morphism $\rho: {\mathfrak g} \to
\End(V)$ of Lie algebras.
\end{defi}

\noindent The functor
$$ \begin{array}{ccccc}
  T_e: & {\mathbf{LieGroup}} & \longrightarrow & {\mathbf{LieAlg}}& \\
   & G & \longmapsto & T_{e}(G)& \\
   & \varphi:G\rightarrow H& \longmapsto & d_{e}\varphi:& T_{e}(G)\rightarrow T_{e}(H)\\
\end{array}$$
induces an equivalence between $S(\mathbf{LieGroup})$ and
$S(\mathbf{LieAlg})$.

\begin{defi}
For any Lie algebra $\mathfrak g$ the {\em adjoint} representation
$\ad: {\mathfrak g} \to \End({\mathfrak g})$ is given by
$\ad(x)(y)=[x,y]$ for all $x,y \in {\mathfrak g}$.
\end{defi}

\begin{defi}

\begin{enumerate}
\item{A subspace $I$ of a Lie algebra ${\mathfrak g}$ is called a
{\em Lie subalgebra} if $[x,y] \in I$ for all $x,y \in I$.}
\item{A subalgebra $I$ of ${\mathfrak g}$ is said to be {\em
abelian} if $[x,y]=0$ for all $x,y \in I$.} \item{A subalgebra $I$
of a Lie algebra ${\mathfrak g}$ is called an {\em ideal} if
$[x,y] \in I$ for all $x \in I$ and $y \in \mathfrak g$.}

\end{enumerate}

\end{defi}

\begin{exa}
For any $k$-algebra the space of derivations of A \\
$\Der(A)=\{d:A
\longrightarrow A \ \ |\ \  d(xy)= d(x)y+xd(y) \mbox{ for all }
x,y \in A \}$is a Lie subalgebra of $\End(A)$.
\end{exa}

\begin{defi}
\begin{enumerate}
\item{A Lie algebra ${\mathfrak g}$ is called {\em simple} if it
has no ideals other than ${\mathfrak g}$ and $\{ 0\}$.} \item{A
Lie algebra ${\mathfrak g}$ is called {\em semisimple} if it has
no abelian ideals other than $\{0\}$.} \item{A maximal abelian
subalgebra ${\mathfrak h}$ of ${\mathfrak g}$ is called a {\em
Cartan subalgebra}.}
\end{enumerate}

\end{defi}
\noindent Next theorem is due to Cartan. A proof of it may be
found in \cite{Fu}.
\begin{thm} Let $\mathfrak g$ be a finite dimensional simple
Lie algebra over $\mathbb{C}$, then $\mathfrak g$ is isomorphic to
one of those in the list $\mathfrak{sl}_{n}(\mathbb{C})$,
$\mathfrak{sp}_{2n}(\mathbb{C})$, $\mathfrak{so}_{2n}(\mathbb{C})$,
$\mathfrak{so}_{2n+1}(\mathbb{C})$, $E_6, E_7, E_8, F_4$ and $G_2$.
\end{thm}
\noindent Algebras $\mathfrak{sl}_{n}(\mathbb{C})$,
$\mathfrak{sp}_{2n}(\mathbb{C})$, $\mathfrak{so}_{2n}(\mathbb{C})$
and $\mathfrak{so}_{2n+1}(\mathbb{C})$ are called classical and will
be explained using our graphical notation in Sections \ref{sl},
\ref{sp}, \ref{sop} and \ref{soi}. The algebras $E_6, E_7, E_8, F_4$
and $G_2$ are called exceptional and the reader may find their
definitions in \cite{Fu}.
\begin{defi}
The {\em Killing form} on $\mathfrak g$ is the bilinear map
$\langle \mbox{ }, \mbox{ } \rangle : {\mathfrak g}\times
{\mathfrak g} \longrightarrow  \mathbb{C}$ given for all $x,y \in
{\mathfrak g}$ by $\langle x,y \rangle= \tr(\ad(x) \circ \ad(y))$,
where $\circ$ denotes the product in $\End({\mathfrak g})$ and
$\tr:\End({\mathfrak g}) \to \mathbb{C}$ is the trace map.
\end{defi}
\noindent Denote by ${\mathfrak h}^{\ast}$ the linear dual of vector
space ${\mathfrak h}$. The following proposition describes
representations of abelian Lie algebras.

\begin{prop}\label{DC}
Let $\mathfrak h$ be an abelian Lie algebra and $\rho: {\mathfrak h}
\to \End(V)$ a representation of $\mathfrak h$. Then $V$
admits a decomposition
\begin{equation}\label{card}
V = \bigoplus_{\alpha \in \Phi} V_{\alpha}\end{equation} where for
each $\alpha \in \mathfrak h^{*}$,  $V_{\alpha}=\{x \in V:
\rho(h)(x)=\alpha(h)x,\ \mbox{ for all }\ \ h\in {\mathfrak h}
\}$, and
 $\Phi=\{ \alpha \in \mathfrak h^{*} |\mbox{  }  \mathfrak
h_{\alpha} \neq 0 \}.$
\end{prop}
\noindent Equation (\ref{card}) is called  {\em Cartan
decomposition} of the representation $\rho$ of $\mathfrak h$.
Proposition \ref{DC} yields a map from $S(\mathbf{LieAlg})$ into
$S(\mathbf{Root})$, which turns out to be a bijection, as follows.
Let ${\mathfrak g}$ be a finite dimensional simple Lie algebra over
$\mathbb{C}$ and $\mathfrak h\subset \mathfrak g$ a Cartan
subalgebra. The killing form $\langle \mbox{ }, \mbox{ }
\rangle:\mathfrak h\times \mathfrak h\longrightarrow\mathbb{C}$
restricted to $\mathfrak h$ is non-degenerated and makes the pair
$(\mathfrak h, \langle \mbox{ } , \mbox{ } \rangle)$ an Euclidean
space. The linear dual $\mathfrak h^{\ast}$ has an induced Euclidean
structure, which we still denote by $\langle \mbox{ } , \mbox{ }
\rangle $ induced by the linear isomorphism $f:\mathfrak h
\longrightarrow \mathfrak h^{\ast}$ given by $f(x)(y)=\langle x,y
\rangle$, for all $x,y\in \mathfrak h$.

\noindent The adjoint representation $\ad:\mathfrak h \longrightarrow
\End(\mathfrak g)$ restricted to $\mathfrak h$ give us a Cartan
decomposition ${\displaystyle \mathfrak g=\bigoplus_{\alpha \in
\Phi}
\mathfrak g_{\alpha}}$ where for each $\alpha\in \mathfrak h$\ \
$\mathfrak g_{\alpha}=\{x\in \mathfrak g|\ \ [h,x]=\alpha(h)x,\
\mbox{ for all}\ \ h\in \mathfrak h \}$ and $\Phi=\{\alpha\in
\mathfrak h^{\ast}: \ \mathfrak h_{\alpha}\neq 0\}$. It is not
difficult to show that the triple $(\mathfrak h^{\ast}, \langle
\mbox{ }, \mbox{ } \rangle,\Phi)$ is a root system.
\begin{defi} Let $\mathfrak g$ be a simple Lie algebra and let
$\Phi$ the root systems associated to $\mathfrak g$. The group $W$
generated by all reflections $S_{\alpha}$ with $\alpha\in\Phi$,
where $S_{\alpha}(\beta)=\beta -2\frac{\langle \alpha , \beta
\rangle} {\langle \alpha , \alpha \rangle}\alpha,$ is known as
the {\em Weyl group} associated to $\mathfrak g$.
\end{defi}
 \noindent One can show that there exists a subset
$\Pi=\{\alpha_1,\alpha_2,\dots,\alpha_n\}$ of $\Phi$ such that
 $\Pi$ is a basis of ${\mathfrak h }^{\ast}$ and
each root $\alpha \in \Phi$  can be written as a linear combination
of roots  in  $\Pi$ with coefficients in $\mathbb{Z}$ which are
either all non-negative or all non-positive. The set $\Pi$ is called
a set of {\em fundamental roots}. The integers
\begin{equation}\label{Cartm}
A_{ij}=2\displaystyle{\frac{\langle \alpha_{i},
\alpha_{j}\rangle}{\langle \alpha_{i}, \alpha_{i}\rangle}}
\end{equation}
are called the {\em Cartan integers} and the matrix $A=(A_{ij})$ is
called the {\em Cartan matrix}. Notice that $A_{ii}=2$, and that for
any $\alpha_{i}, \alpha_{j}\in \Pi$ with $i\neq j$,
$S_{\alpha_{i}}(\alpha_{j})$ is a $\mathbb{Z}$-combination of
$\alpha_{i}$ and $\alpha_{j}$. Since the coefficient of $\alpha_{j}$
is $1$, the coefficient associated to $\alpha_{i}$ in
$S_{\alpha_{i}}(\alpha_{j})$  must be a non-positive integer, i.e.,
$A_{ij}\in \mathbb{Z}^{\leq 0}$. The angle $\theta_{ij}$ between
$\alpha_{i}, \ \alpha_{j}$ is given by the cosine formula
$$\langle \alpha_{i},
\alpha_{j}\rangle=\langle \alpha_{i},\alpha_{i} \rangle
^{\frac{1}{2}} \langle \alpha_{j},
\alpha_{j}\rangle^{\frac{1}{2}}\cos(\theta_{ij}).$$  Then we have
$$4 \cos^{2}(\theta_{ij})=2\frac{\langle\alpha_{i},\alpha_{j}\rangle}
{\langle\alpha_{j},\alpha_{j}\rangle}\cdot
2\frac{\langle\alpha_{j},\alpha_{i}\rangle}{\langle\alpha_{i},\alpha_{i}\rangle},$$
and therefore $4 \cos^{2}(\theta_{ij})=A_{ij}A_{ji}.$ Let
$n_{ij}=A_{ij}A_{ji}$ clearly $n_{ij}\in \mathbb{Z}$ and $n_{ij}\geq
0$. Since $-1\leq \cos(\theta_{ij})\leq 1$ the only possible values
for $n_{ij}$ are $n_{ij}=0,1,2\ \mbox{or}\ 3$.

\begin{defi}\label{DD}
The Dynkin diagram $\Delta$ associated to a simple Lie algebra
${\mathfrak g}$ is the graph $\Delta$ with vertices
$\{1,\dots,n\}$ in bijective correspondence  with the set $\Pi$ of
fundamental roots of ${\mathfrak g}$ such that
\begin{enumerate}
\item{Vertices $i,j$ with $i\neq j$ are joined by $n_{ij}=A_{ij}A_{i}$ edges, where
$A_{ij}$ is given by formula $(\ref{Cartm})$.}
\item{Between each double edge or triple we attach the symbol\  $<$,
 or the symbol\   $>$  \ pointing towards the shorter root with
respect to Killing form.}
\end{enumerate}
\end{defi}
\begin{thm}\label{iso}
Consider the root decomposition of the simple Lie algebra $\mathfrak
g$ with respect to $\mathfrak h$. Let $\alpha\in \Phi$ be a root,
for each nonzero $x_{\alpha}\in \mathfrak g_{\alpha}$ there is
$x_{-\alpha}\in \mathfrak g_{\alpha}$ and $ h_{\alpha}\in \mathfrak
h$ such that $\alpha(h_{\alpha})=2$, $[x_{\alpha},x_{-\alpha}]=
h_{\alpha}$, $[ h_{\alpha},x_{\alpha} ]=2x_{\alpha}$ and $[
h_{\alpha},x_{-\alpha} ]=2x_{-\alpha}$.
\end{thm}

\noindent Theorem \ref{iso} implies that for any root
$\alpha\in\Phi$, $x_{\alpha}, x_{-\alpha}$ and $h_{\alpha}$ span a
subalgebra $\mathfrak s_{\alpha}$ such that $\mathfrak
s_{\alpha}\cong \mathfrak {sl}_{2}(\mathbb{C})$, see Section
\ref{sl2} for a definition of $\mathfrak {sl}_{2}(\mathbb{C})$. This
fact explain the distinguished role played by  $\mathfrak
{sl}_{2}(\mathbb{C})$ in the representation theory of simple Lie
algebras.
\begin{defi}
Let $\Phi\subset\mathfrak{h}$ be the root system associated with a
simple Lie algebra $\mathfrak{g}$.

\begin{enumerate}
\item{For any $\alpha \in \Phi$, the Cartan element $h_{\alpha}\in \mathfrak{h}$  given by
Theorem \ref{iso} is called the  {\em coroot} associated to root
$\alpha$. $\Phi_c=\{h_{\alpha}: \alpha\in \Phi\}$ is the coroot
system associated to ${\mathfrak h}$ and $\Pi_c=\{
h_{\alpha}:\alpha\in\Pi\}$ is the set of fundamental coroots. }

\item{The elements
$w_1,\dots,w_n$ in $\mathfrak{h}^{\ast}$ given by the relations
$w_{i}(h_{j})=\delta_{ij}$, for all $1\leq i,j\leq n$, where $h_{j}$
is the coroot associated to fundamental root $\alpha_{j}$, are
called the fundamental weights.}
\end{enumerate}
\end{defi}

\noindent One can recover a simple Lie algebra  $\mathfrak g$ from its associated Dynkin
diagram $\Delta$ as follows: Let $n_{ij}$ be the adjacency matrix of
$\Delta$. The relation $n_{ij}=A_{ij}A_{ji}$ determines univocally
the Cartan matrix $A_{ij}$. Consider the free Lie algebra generated
by the symbols $H_1,\dots,H_n, X_1,\dots,X_n,Y_1,\dots,Y_n$. Form
the quotient of this free Lie algebra  by the relations
$$
\begin{array}{llc}
  [H_i,H_j]=0\ \ (\mbox{all} \ i,j); & [X_i,Y_i]=H_i \ \ (\mbox{all} \ i);
  & [X_i,X_j]=0 \ \ ( i\neq j); \\
  \mbox{} & \mbox{}&\mbox{}\\

  [H_i,X_j]=A_{ij}X_j\ \ (\mbox{all} \ i,j); & [H_i,Y_j]=-A_{ij}Y_{j}\ \ (\mbox{all} \ i,j); &  \mbox{ }\\
\end{array}
$$

and for all $i\neq j$,
$$\begin{array}{lll}
  [X_i,X_j]=0, & [Y_i,Y_j]=0, & \mbox{if}\ \ A_{ij}=0. \\
  \mbox{} & \mbox{} & \mbox{} \\

  [X_i,[X_i,X_j]]=0, & [Y_i,[Y_i,Y_j]]=0 & \mbox{if}\ \ A_{ij}=-1. \\
  \mbox{} & \mbox{} & \mbox{} \\

  [X_i,[X_i,[X_i,X_j]]]=0, & [Y_i,[Y_i,[Y_i,Y_j]]]=0 & \mbox{if}\ \ A_{ij}=-2. \\
  \mbox{} & \mbox{} & \mbox{} \\

  [X_i,[X_i,[X_i,[X_i,X_j]]]]=0, & [Y_i,[Y_i,[Y_i,[Y_i,Y_j]]]]=0 & \mbox{if}\ \ A_{ij}=-3. \\
\end{array}
$$

\noindent Serre shows that the resulting Lie algebra is a
finite-dimensional simple Lie algebra isomorphic to $\mathfrak g$.
See
\cite{SE} for more details.
\subsection {Jacobian Criterion}

Let $k$ be a field of characteristic zero and let $V$ be a finite
dimensional $k$-vector space. Set $V=\langle e_{1},e_{2},\dots,
e_{n}\rangle$ and $V^{\ast}=\langle x_{1},x_{2}, \dots ,
x_{n}\rangle$ such that $x_{i}(e_{j})=\delta_{ij}$. We have
$$x_{k}\Big(  \sum_{i=1}^{n} a_{i} e_{i} \Big)= a_{k}, \ \ 1\leq k \leq n.$$
Any $f\in V^{\ast}$ is written as  $f= b_{1}x_{1}+b_{2}x_{2}+\dots
+b_{n}x_{n}$. Denote by $S=S(V^{\ast})$ the symmetric algebra of the
dual space $V^{\ast}$ which can be identify with the polynomial ring
$k[x_{1},x_{2},\dots, x_{n}]$. Let $G$ be a finite group which acts
on $V$. $G$ also acts on $V^{\ast}$, and thus it acts on
$S=S(V^{\ast})$ as follows
$$\begin{array}{ccc}
  S(V^{\ast})\times G & \longrightarrow & S(V^{\ast}) \\
  (p,g) & \longmapsto & p(g)
\end{array}$$
where $(p g)(v)=p(gv),\ \ \mbox{for all}\ \ g\in G,
\ p\in S(V^{\ast}),\  v\in V.$ The algebra
$$k[x_{1},x_{2},\dots,x_{n}]^{G}=\{ p\in k[x_{1},x_{2},\dots,x_{n}]: \  p(g)=p,\
 \forall g\in G\}$$
is called the $G$-invariant subalgebra of
$k[x_{1},x_{2},\dots,x_{n}]$.
\begin{defi}
Let $K$ be a field and $F$ a extension of $K$. Let $S$ be a subset
of  $F$. The set $S$ is algebraically dependent over $K$ if for
any positive integer $n$ there is a polynomial non-zero $f\in
k[x_1,\dots,x_n]$ such that $f(s_1,\dots,s_n)=0$ for any different
$s_1,\dots,s_n\in S$. In otherwise $S$ is algebraically
independent.
\end{defi}
\begin{thm}\label{Winv}
Let $\mathbb{C}[x_{1},x_{2},\dots,x_{n}]^{W}$ be the subalgebra of
$\mathbb{C}[x_{1},x_{2},\dots,x_{n}]$ consisting of $W$-invariant
polynomials, then $\mathbb{C}[x_{1},x_{2},\dots,x_{n}]^{W}$ is
generated as an  $\mathbb{C}$-algebra by $n$ homogeneous,
algebraically independent elements of positive degree together
with $1$.
\end{thm}
\noindent The idea of proof of Theorem \ref{Winv} goes as follows:
let $I$ be the ideal of $\mathbb{C}[x_{1},x_{2},\dots,x_{n}]$
generated by all homogeneous $W$-invariant polynomials of positive
degree. Using Hilbert's Basis Theorem we may choose a minimal
generating set $f_{1},f_{2},\dots, f_{r}$ for $I$  consisting of
homogeneous $W$-invariant polynomials of positive degree. One can
show that $r=n$ and furthermore
$\mathbb{C}[x_{1},x_{2},\dots,x_{n}]^{W}=
\mathbb{C}[f_{1},f_{2},\dots,f_{n}]$.

\begin{prop}
Let $f_1,\dots, f_n$ and  $g_1,\dots, g_n$ be two sets of
homogeneous, algebraically independent generators  of
$\mathbb{C}[x_1,\dots,x_n]^{W}$  with degrees $d_i$ and $e_i$
respectively, then (after reordering) $d_i=e_i$ for all
$i=1,\dots,n$.
\end{prop}
\noindent The numbers $d_1,\dots,d_n$ written in increasing order
are called the {\em degrees} of  $W$.  Theorem \ref{CJ} below is a
simple criterion for the algebraic independence of polynomials
$f_1,\dots,f_n$ expressed in terms of the Jacobian determinant. We
write $J(f_1,\dots,f_n)$ for the determinant of the $n\times n$
matrix whose $(i,j)$-entry is $\frac{\partial f_i}{\partial x_j}$.

\begin{thm}[Jacobian criterion]\label{CJ}
The set of polynomials $f_1,\dots,f_n\in k[x_1,\dots,x_n]$ are
algebraically independent over a field $k$ of characteristic zero
if and only if $J(f_1,\dots,f_n)\neq 0$.
\end{thm}

\section{Graph and matrices}

We denote by $\textbf{Digraph}^{1}(n,n)$ the vector space
generated by bipartite directed graphs with a unique edge starting
on the set $[n]$ and ending on the set $[n]$.  We describe
$\textbf{Digraph}^{1}(n,n)$ pictorially as follows
\begin{figure}[ht]
\begin{center}
\includegraphics[width=3in]{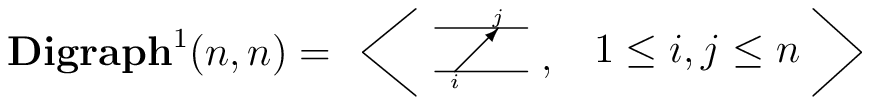}
\end{center}
\end{figure}

\noindent where the symbol
\begin{figure}[ht]
\begin{center}
\includegraphics[width=0.4in]{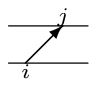}
\end{center}
\end{figure}

\noindent denotes the graph whose unique edge stars at vertex at $i$
and ends at vertex $j$. We define a product on
$\textbf{Digraph}^{1}(n,n)$ as follows
\begin{figure}[ht]
\begin{center}
\includegraphics[width=3.5in]{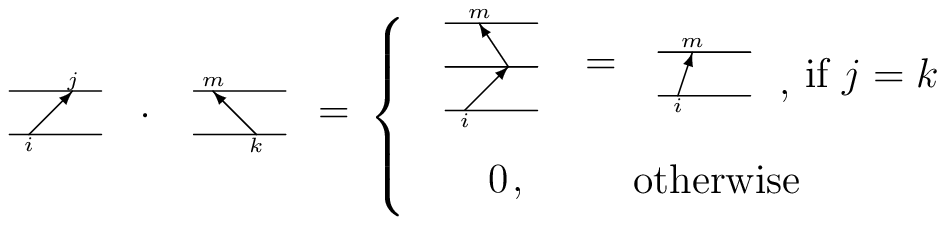}
\end{center}
\end{figure}

\noindent The trace on $\textbf{Digraph}^{1}(n,n)$ is defined as the
linear functional $\tr: \textbf{Digraph}^{1}(n,n) \longrightarrow
\mathbb{C}$ given by

\begin{figure}[!h]
\begin{center}
\includegraphics[width=3in]{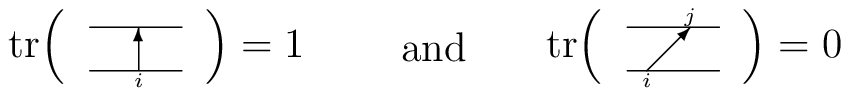}
\end{center}
\end{figure}

\noindent Algebra $\textbf{Digraph}^{1}(n,n)$ is isomorphic to
$\End(\mathbb{C}^{n})$ through the application
\begin{figure}[!h]
\begin{center}
\includegraphics[width=2.4in]{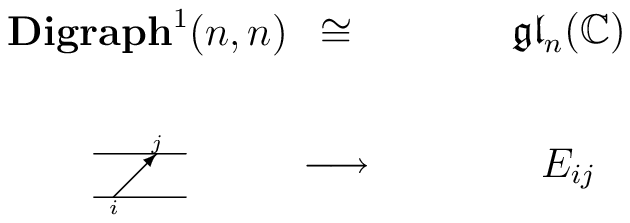}
\end{center}
\end{figure}

\noindent We will use this isomorphism to give combinatorial
interpretation of results on $\End(\mathbb{C}^{n})$ that are
traditionally express in the language of matrices.

\section{Linear special algebra ${\mathfrak
{sl}}_{2}(\mathbb{C})$}\label{sl2}

We begin studying the special linear algebra ${\mathfrak
{sl}}_{2}(\mathbb{C})$. It plays a distinguished role in the theory
of Lie algebras. By definition we have
$${\mathfrak {sl}}_2(\mathbb{C})=\{A\in \End(\mathbb{C}^{2}): \tr(A)=0\}.$$
As subspace of $\textbf{Digraph}^{1}(2,2)$, ${\mathfrak
{sl}}_{2}(\mathbb{C})$ is the following vector space
\begin{figure}[!h]
\begin{center}
\includegraphics[width=3.3in]{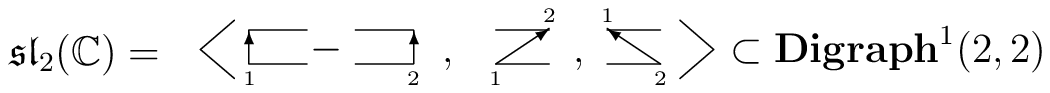}
\end{center}
\end{figure}

\noindent We fix as Cartan subalgebra of ${\mathfrak
{sl}_2(\mathbb{C})}$ the $1$-dimensional subspace

\begin{figure}[!h]
\begin{center}
\includegraphics[width=2.7in]{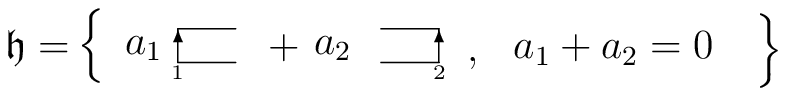}
\end{center}
\end{figure}

\noindent The dual space is ${\mathfrak h}^{\ast}=\langle
a_1,a_2\rangle / \{a_1+a_2=0\}$, where

\begin{figure}[!h]
\begin{center}
\includegraphics[width=1in]{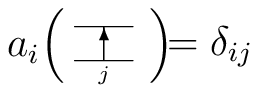}
\end{center}
\end{figure}

\subsection{Root system of ${\mathfrak {sl}}_2(\mathbb{C})$}

Consider the projection map $\langle a_1,a_{2}\rangle
\longrightarrow \langle a_1, a_2\rangle /\{a_1+a_2=0\}$. We still
denote by $a_i$ the image of $a_i$ under the projection above. Now
each $h\in{\mathfrak{h}}$ is of the form
\begin{figure}[!h]
\begin{center}
\includegraphics[width=2.5in]{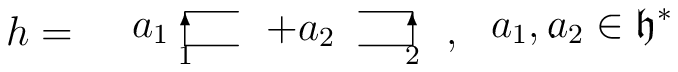}
\end{center}
\end{figure}
\vspace{0.2cm}

\noindent Let us compute the roots


\begin{figure}[!h]
\begin{center}
\includegraphics[width=2in]{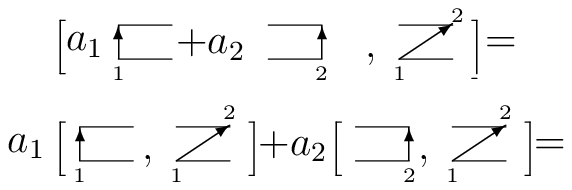}
\includegraphics[width=3.5in]{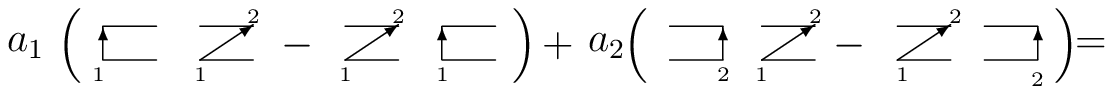}
\includegraphics[width=2in]{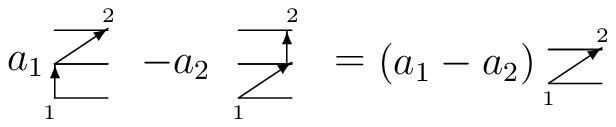}
\end{center}
\end{figure}




\noindent Therefore the root system of ${\mathfrak
{sl}}_2(\mathbb{C})$ is $\Phi=\{a_1-a_2,a_2-a_1\}.$ Setting
$\alpha=a_1-a_2$ we have that the roots are $\alpha$ and $-\alpha$
and the set of fundamental roots is $\Pi=\{\alpha\}.$ In pictures
\begin{figure}[ht]
\begin{center}
\includegraphics[width=1.5cm]{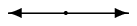}
\end{center}
\end{figure}

\subsection{Coroot system of ${\mathfrak {sl}_2(\mathbb{C})}$}

Let $x_{\alpha}$ and $x_{-{\alpha}}$ be the covectors associated
with the roots $\alpha$ and $-\alpha$ of ${\mathfrak
{sl}}_2(\mathbb{C})$ respectively.
\begin{figure}[!h]
\begin{center}
\includegraphics[width=10cm]{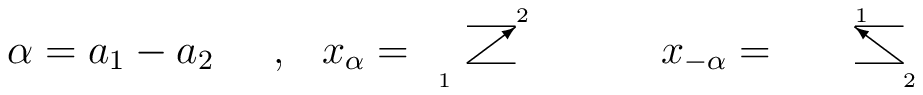}
\end{center}
\end{figure}

\noindent A vector $h_{\alpha}\in {\mathfrak h}$ is said to be the
coroot associated to the root $\alpha\in {\mathfrak h}^{\ast}$, if
$h_{\alpha}=c[x_{\alpha}, x_{-\alpha}]$, $c\in \mathbb{C}$ and
$\alpha(h_{\alpha})=2$.

\begin{figure}[!h]
\begin{center}
\includegraphics[width=8cm]{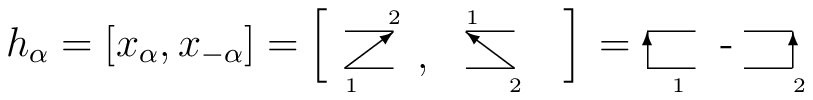}
\end{center}
\end{figure}
\noindent since
\begin{figure}[!h]
\begin{center}
\includegraphics[width=2in]{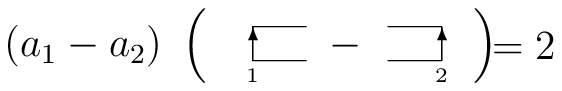}
\end{center}
\end{figure}

\subsection{Killing form of ${\mathfrak {sl}}_2(\mathbb{C})$}
\begin{eqnarray*}
\langle x,y\rangle & = & \sum_{\alpha\in \Phi}\alpha(x)\alpha(y)\\
\mbox{ } & = &2(x_1-x_2)(y_1-y_2)\\
\mbox{ } & = & 2x_1y_1+2x_2y_2-2x_1y_2-2x_2y_1\\
\mbox{ } & = &
2(x_1y_1+x_2y_2)-2(x_1+x_2)(y_1+y_2)+2(x_1y_1+x_2y_2)\\
\mbox{ } & = &4\tr(xy).
\end{eqnarray*}

\subsection{Dynkin diagram of ${\mathfrak
{sl}}_2(\mathbb{C})$}

We have only one fundamental root, so the Dynkin diagram is just
$\circ$.

\section{Special linear algebra ${\mathfrak {sl}}_{n}(\mathbb{C})$}\label{sl}

Let us recall the special linear algebra ${\mathfrak
{sl}}_{n}(\mathbb{C})$
$${\mathfrak {sl}}_{n}(\mathbb{C})=\{ A\in \End(\mathbb{C}^{n}):\ \tr(A)=0\}$$
${\mathfrak {sl}}_{n}( \mathbb{C} )$ consider as a subspace of
$\textbf{Digraph}^{1}(n,n)$ is following subspace

\begin{figure}[!h]
\begin{center}
\includegraphics[width=4.3in]{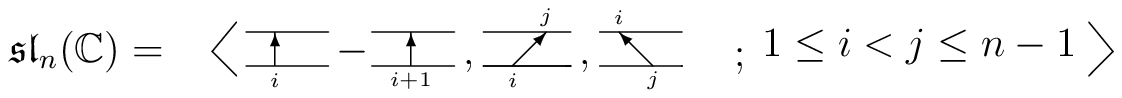}
\end{center}
\end{figure}

\subsection{Root system of ${\mathfrak {sl}}_{n}(\mathbb{C})$}

We take as Cartan subalgebra the subspace of ${\mathfrak
{sl}}_{n}(\mathbb{C})$
\vspace{-1cm}
\begin{figure}[!h]
\begin{center}
\includegraphics[width=4.5in]{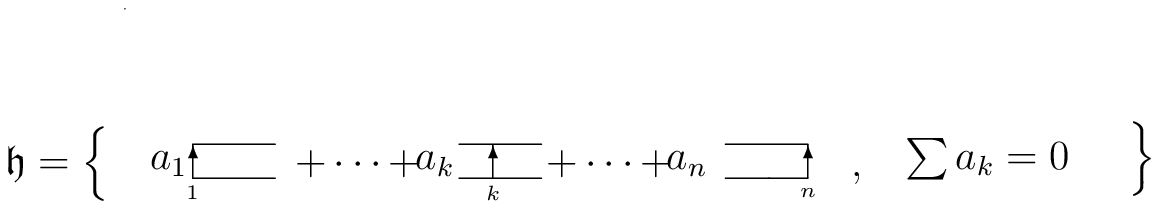}
\end{center}
\end{figure}

\noindent The dual space is ${\mathfrak h}^{\ast}=\langle
a_1,\dots,a_{n}\rangle/(\sum a_i=0)$, where

\begin{figure}[!h]
\begin{center}
\includegraphics[width=1in]{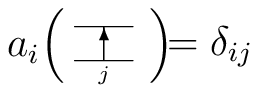}
\end{center}
\end{figure}
\noindent Consider the projection $\langle a_1,\dots,a_{n}\rangle
\longrightarrow \langle a_1,\dots,a_{n}\rangle /(\sum a_k=0)$. The
image of $a_i$ under the projection above is still denote by $a_i$.
Then vector $h \in {\mathfrak h}$ can be written as
\vspace{-1cm}
\begin{figure}[!h]
\begin{center}
\includegraphics[width=2.8in]{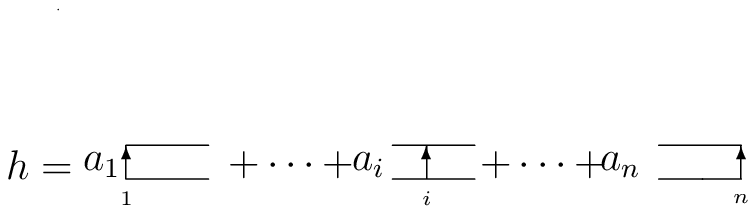}
\end{center}
\end{figure}

\noindent Let us compute the root system
\vspace{-0.5cm}
\begin{figure}[!h]
\begin{center}
\includegraphics[width=3in]{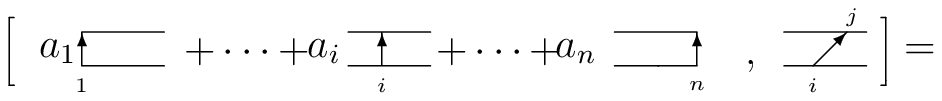}
\includegraphics[width=2.5in]{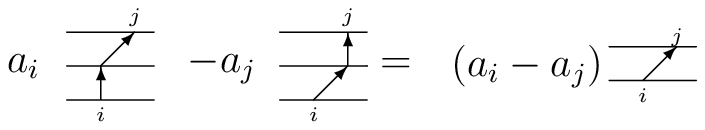}
\end{center}
\end{figure}

\noindent Also
\vspace{-0.5cm}
\begin{figure}[!h]
\begin{center}
\includegraphics[width=3in]{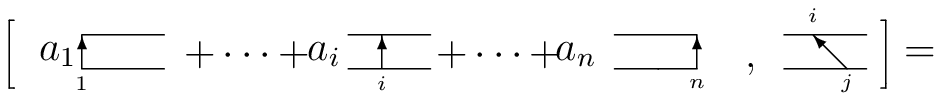}
\includegraphics[width=2.5in]{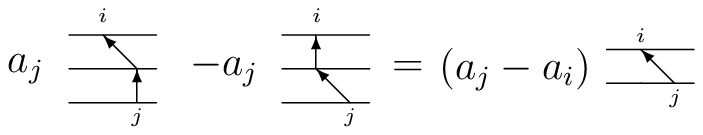}
\end{center}
\end{figure}

\noindent Thus the root system of ${\mathfrak{sl}}_{n}(\mathbb{C})$ is
$\Phi=\{a_i-a_j,\ a_j-a_i,\ 1\leq i<j\leq n-1\}\subset {\mathfrak
h}^{\ast}$. The set of fundamental roots is $\Pi=\{a_i-a_{i+1},
\ i=1,\dots, n-1\}$. In pictures for $n=2,3\ \mbox{and}\ 4$ the root
systems look like
\begin{figure}[!h]
\begin{center}
\includegraphics[width=2.8in]{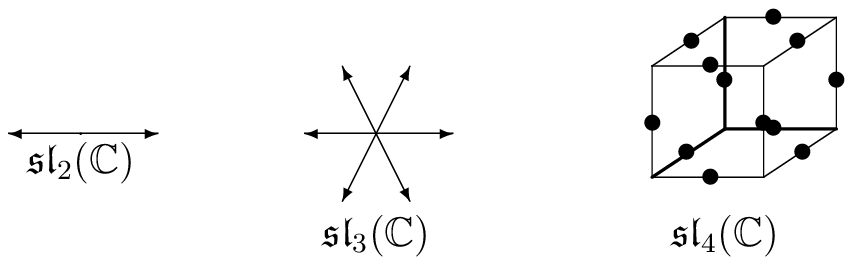}
\end{center}
\end{figure}
\newpage

\noindent Consider  the linear map
$$T: \textbf{Digraph}^{1}(n,n) \longrightarrow
\textbf{Digraph}^{1}(n,n)$$ sending each directed graph into its
opposite graph. Clearly $T$ is an antimorphism, i.e,
$T(ab)=T(b)T(a)$, for all $a,b\in \textbf{Digraph}^{1}(n,n)$. For
example,

\begin{figure}[ht]
\begin{center}
\includegraphics[width=1.5in]{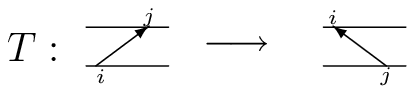}
\end{center}
\end{figure}

\noindent Notice that negative roots can be obtain from the positive
ones through an application of $T$.

\subsection{Coroots and weights for ${\mathfrak {sl}}_{n}(\mathbb{C})$}

\begin{enumerate}
\item Coroot associated to the root $a_i-a_j$

\begin{figure}[!h]
\begin{center}
\includegraphics[width=2in]{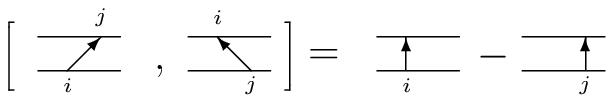}
\end{center}
\end{figure}

\item Coroot associated to the root $a_j-a_i$

\begin{figure}[!h]
\begin{center}
\includegraphics[width=2.5in]{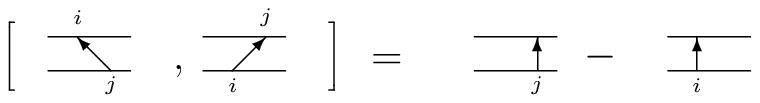}
\end{center}
\end{figure}

\end{enumerate}

\noindent The set of fundamental coroots has the form
$\Pi_{c}=\{h_{i}-h_{i+1}, 1\leq i\leq n-1\}$ where

\begin{figure}[!h]
\begin{center}
\includegraphics[width=2in]{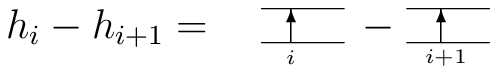}
\end{center}
\end{figure}

\noindent The set of fundamental weights is $w_{i}=a_1
+a_2+\cdots+a_i$ since

\begin{figure}[!h]
\begin{center}
\includegraphics[width=3in]{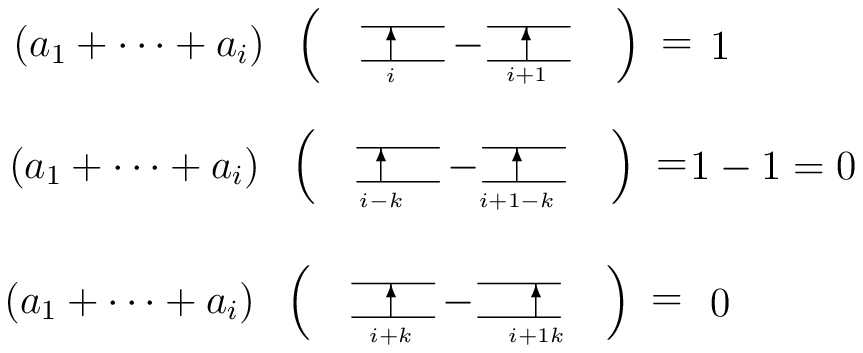}
\end{center}
\end{figure}
 \newpage

\subsection{The Killing form of ${\mathfrak {sl}}_{n}(\mathbb{C})$}

Let $x$ and $y$ in ${\mathfrak h}$. Set

\vspace{-0.5cm}
\begin{eqnarray*}
\langle x,y\rangle & = & \sum_{\alpha\in \Phi}\alpha(x)\alpha(y)\\
\mbox{ }& = &\sum_{i<j}(x_i-x_j)(y_i-y_j)+\sum_{i<j}(x_j-x_i)(y_j-y_i)\\
\mbox{ }& = &2\sum_{i<j}(x_i-x_j)(y_i-y_j)\\
\mbox{ }& = &2 \big( \sum_{i<j}x_iy_i+\sum_{i<j}x_jy_j-\sum_{i<j}x_iy_j-\sum_{i<j}x_jy_i \big)\\
\mbox{ }& = &2 \big( \sum (n-i)x_iy_i+ \sum (i-1)x_iy_i +\sum x_iy_i \big)\\
\mbox{ }& = &2n\tr(xy).
\end{eqnarray*}

\subsection{Weyl group of ${\mathfrak {sl}}_{n}(\mathbb{C})$}

Consider the fundamental roots $\alpha_i=a_i-a_{i+1}, \ \
i=1,\dots,n-1$ and $S_{\alpha_{i}}$ the reflection associated to the
fundamental root ${\alpha_{i}}$. Let $h\in{\mathfrak h}$ and
$h_{\alpha_i}$  be the coroot associated to the fundamental root
$\alpha_i$. By definition we have $S_{\alpha_{i}}
(h)=h-\alpha_{i}(h)h_{\alpha_i}$

\begin{figure}[!h]
\begin{center}
\includegraphics[width=5.6in]{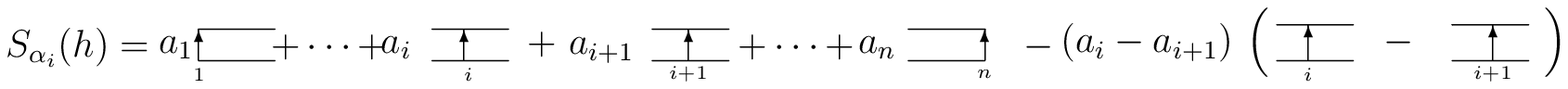}
\vspace{0.5cm}
\includegraphics[width=3.6in]{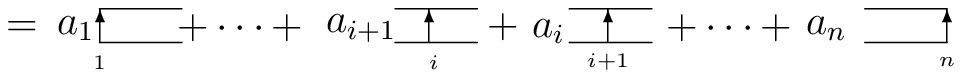}
\end{center}
\end{figure}

\noindent so we see that reflections $S_{\alpha_{i}}$ has the form

\begin{figure}[!h]
\begin{center}
\includegraphics[width=2.3in]{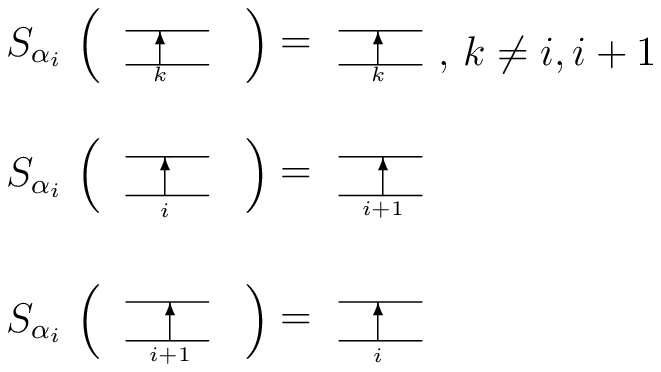}
\end{center}
\end{figure}

\noindent Therefore the Weyl group $A_n$ associated with ${\mathfrak
{sl}}_{n+1}(\mathbb{C})$ is the symmetric group on $n$ letters
$$A_{n}=\langle S_{\alpha_{i}}|\ \ i=1,\dots,n-1\rangle=S_{n}.$$

\subsection{Dynkin diagram of ${\mathfrak {sl}}_{n}(\mathbb{C})$ and Cartan matrix.}

Using equation (\ref{Cartm}) one can checks that the Cartan matrix
associated to the Lie algebra ${\mathfrak {sl}}_{n}(\mathbb{C})$ is

$$A_n= \left( \begin{array}{cccccc}
  2 & \!\!\!\!-1 & 0 & 0 & \ldots& 0 \\
 \!\!\!\! -1 & 2 &\! \!\!\!-1 & 0 &   \ldots& 0\\
    0 & \!\!\!\!-1 & 2 & \!\!\!\!-1 &  \ldots& 0   \\
   &  & & \ddots  &  &    \\
   &  &  &   & \ddots &    \\
  0& 0  & \ldots & 0 &\! \!\!\!-1 & 2
\end{array}\right) $$

\noindent  The Dynkin diagram associated to ${\mathfrak {sl}}_{n}(\mathbb{C})$
is
\begin{figure}[ht]
\begin{center}
\includegraphics[width=2.5in]{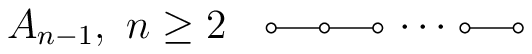}
\end{center}
\end{figure}
\subsection{Invariant polynomials for ${\mathfrak {sl}}_{n}(\mathbb{C})$}

Consider the action of $S_{n+1}$ on $\mathbb{R}^{n+1}$ given by
$$\begin{array}{ccc}
S_{n+1}\times  {\mathbb{R}}^{n+1} & \longrightarrow &{\mathbb{R}^{n+1}}  \\
(\pi,x)& \longmapsto & (\pi x)_i=x_{{\pi}^{-1}(i)}
\end{array}$$
notice that the permutation $(ij)$ acts as a reflection on
$\mathbb{R}^{n+1}$ since
$$(ij)(x_i-x_j)=x_j-x_i=-(x_i-x_j)$$
$$(ij)(x)=x,\ \ \mbox{si}\  \ x\in(x_i-x_j)^{\bot}\ (\mbox{es decir}\  x_i=x_j )$$
Since $S_{n+1}$ is generated by transpositions $(i\ i+1), \ i=1,
\dots,n $, then $S_{n+1}$ is an example of what is called a
reflection group. Recall that a linear action of a group $G$ on a
vector space $V$ is said to be effective if the only fixed point
is $0$. The action of $S_{n+1}$ on $\mathbb{R}^{n+1}$ fixes points
in $\mathbb{R}^{n+1}$ lying on the straight line $\{(x,x,\dots,x)|
\ \ x\in\mathbb{R}\}$. Thus the action of $A_n$ on
$\mathbb{R}^{n+1}$ fails to be effective. If we instead let $A_n$
act on the hyperplane $V=\{(x_1 \dots x_{n+1}) \in
\mathbb{R}^{n+1}|\  x_1+\dots+x_{n+1}=0\}$ then the action becomes
effective. Consider the power symmetric functions
$$f_i=x_1^{i+1}+\cdots+x_{n+1}^{i+1},\ \ 1\leq i\leq n.$$
Each $f_i$ is $S_{n+1}$-invariant, and together the power
symmetric functions form a set of basic invariants. This fact can
be proven as follows: first notice that

$$\gr(f_1) \gr(f_2)\cdots \gr(f_n)=2\cdot3 \cdots n(n+1)=(n+1)!=|S_{n+1}|=|A_n|.$$
Next, it is easy to compute the Jacobian $J(f_1,f_2, \cdots ,f_n)$
yielding the nonvanishing polynomial
$$J(f_1,f_2, \cdots ,f_n)={\displaystyle (n+1)!\prod_{1\leq i<j\leq n}(x_j-x_i)\prod_{i=1}^{n}
(x_1+\cdots+2x_i+\cdots +x_n)}$$ Finally, use the Jacobian
criterion.

\section{Symplectic Lie algebra ${\mathfrak
{sp}_{2n}}(\mathbb{C})$}\label{sp}

Recall that the symplectic Lie algebra $\mathfrak {sp}_{2n}$ is
defined as

$${\mathfrak {sp}}_{2n}(\mathbb{C})=\{X:X^{t}S+SX=0\}.$$ Here  $S\in
M_{2n}(\mathbb{C})$ is the matrix
$$S=\left( \begin{array}{cc}
  0 & I_n \\
  -I_n & 0
\end{array}\right) $$
Equivalently,
$${\mathfrak {sp}}_{2n}(\mathbb{C})=\Bigg \{\left(
\begin{array}{cc}
  A & B \\
  C & -A^{t}
\end{array}\right);\ A,B,C\in M_{n}(\mathbb{C})\  \mbox{y}\   B=B^{t},\  C=C^{t} \Bigg \}$$
${\mathfrak {sp}_{2n}}(\mathbb{C})$ as a subspace of
$\textbf{Digraph}^{1}(2n,2n)$ is given by
\begin{figure}[!h]
\begin{center}
\includegraphics[width=6in]{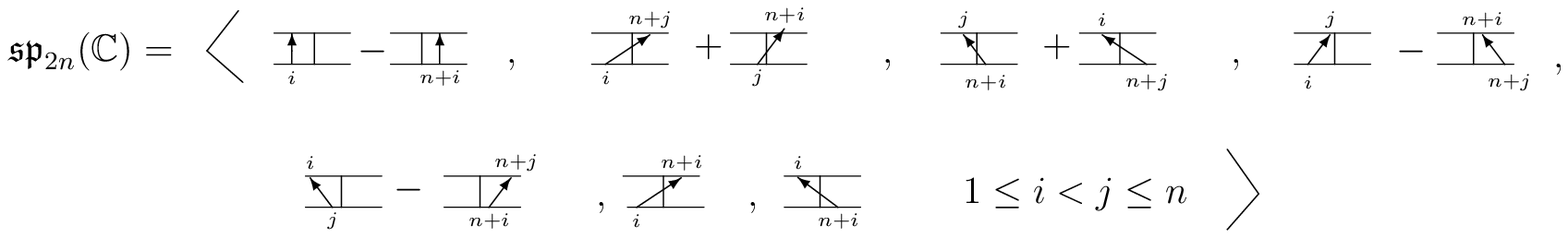}
\end{center}
\end{figure}

\noindent We take as a Cartan subalgebra of ${\mathfrak
{sp}}_{2n}(\mathbb{C})$ the following  subspace

\begin{figure}[!h]
\begin{center}
\includegraphics[width=3in]{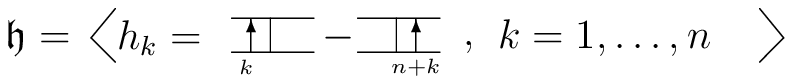}
\end{center}
\end{figure}

\subsection{Root system of ${\mathfrak {sp}}_{2n}(\mathbb{C})$}

Consider $h\in {\mathfrak h}$

\begin{figure}[!h]
\begin{center}
\includegraphics[width=2in]{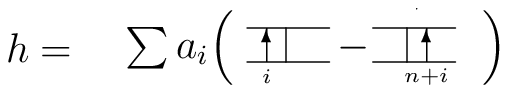}
\end{center}
\end{figure}

\noindent where $ \{a_i\}$ denotes de base of ${\mathfrak h}^{\ast}$
dual to the given base of ${\mathfrak h}$.  Let us define $$T:
\textbf{Digraph}^{1}(2n,2n) \longrightarrow
\textbf{Digraph}^{1}(2n,2n)$$ to be the linear map that sends each
directed graph into its opposite. Clearly $T$ es un antimorphism,
i.e, $T(ab)=T(b)T(a)$ for all $a,b\in \textbf{Digraph}^{1}(2n,2n)$.
For example,
\begin{figure}[ht]
\begin{center}
\includegraphics[width=1.5in]{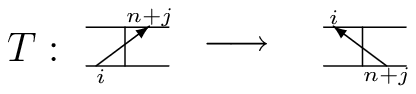}
\end{center}
\end{figure}
We will compute explicitly the positive roots. To obtain the
negative roots it is enough to apply the transformation $T$ to
each positive root.

\begin{figure}[!h]
\begin{center}
\includegraphics[width=3.4in]{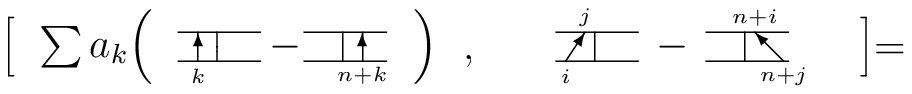}
\includegraphics[width=1.7in]{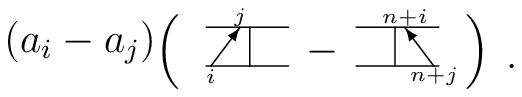}
\includegraphics[width=3.2in]{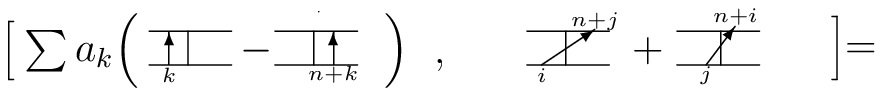}
\includegraphics[width=1.9in]{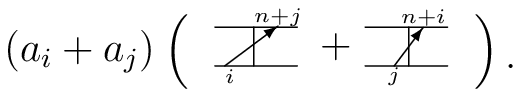}
\includegraphics[width=3.2in]{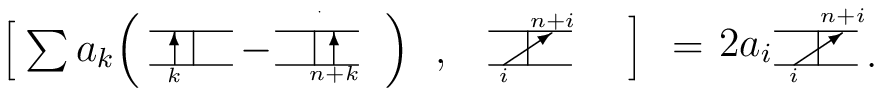}
\end{center}
\end{figure}
\noindent Thus the root system of $ {\mathfrak
{sp}}_{2n}(\mathbb{C})$ is $\Phi=\{ a_i-a_j, a_j-a_i, a_i+a_j,
-a_i-a_j, 2a_i, -2a_i \ \ 1\leq i<j\leq n\}.$ The set of fundamental
roots is $\Pi=\{\alpha_1, \dots,\alpha_n\}$ where
$\alpha_i=a_i-a_{i+1},  i=1,\dots,n-1$ and $\alpha_n=2a_{n}$. In
pictures, the root system of ${\mathfrak {sp}}_{6}(\mathbb{C})$\quad
looks like
\begin{figure}[!h]
\begin{center}
\includegraphics[height=2cm]{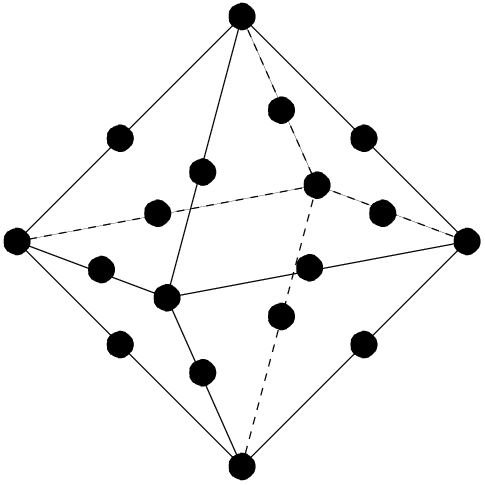}
\end{center}
\end{figure}

\subsection{Coroots and weights of ${\mathfrak {sp}}_{2n}(\mathbb{C})$}

\begin{enumerate}
\item Coroot  asociated to the root $a_i-a_j$
\begin{figure}[!h]
\begin{center}
\includegraphics[width=2.5in]{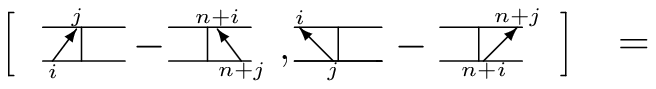}
\includegraphics[width=2.5in]{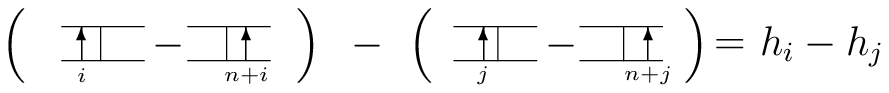}
\end{center}
\end{figure}

$h_i-h_j$ is the coroot associated to the root $a_i-a_j$, since
$(a_i-a_j)(h_i-h_j)=2$.

\item Coroot associated to the root $a_i+a_j$
\newpage
\begin{figure}[!h]
\begin{center}
\includegraphics[width=2.5in]{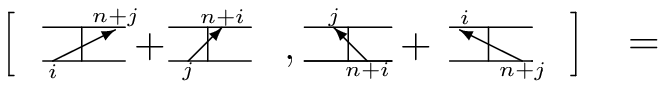}
\includegraphics[width=3.2in]{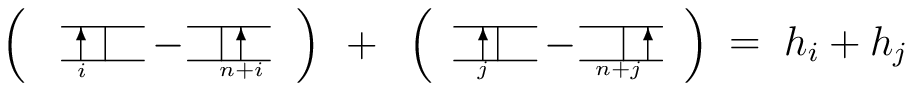}
\end{center}
\end{figure}
$h_i+h_j$ is the coroot associated to the la root $a_i+a_j$, since
$(a_i+a_j)(h_i+h_j)=2$.

\item{Coroot associated to the root $2a_i$

\begin{figure}[!h]
\begin{center}
\includegraphics[width=3.2in]{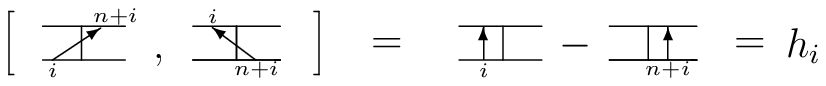}
\end{center}
\end{figure}
}
\end{enumerate}

\noindent We conclude that $\Phi_{c}=\{h_i-h_j, h_j-h_i, h_i+h_j, -h_i-h_j,h_i,-h_i,\
1\leq i<j\leq n\}$ is the coroot system of ${\mathfrak
{sp}}_{2n}(\mathbb{C})$. The set of fundamental coroots is given by
$\Pi_{c}=\{h_i-h_{i+1}, h_n ;\ \  1\leq i\leq n-1\}$ where

\begin{figure}[!h]
\begin{center}
\includegraphics[width=3.3in]{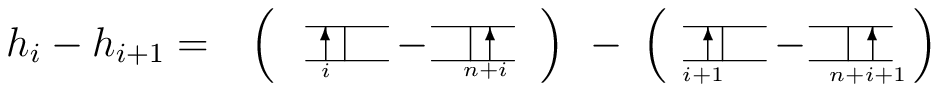}
\end{center}
\end{figure}

\noindent and

\begin{figure}[!h]
\begin{center}
\includegraphics[width=1.3in]{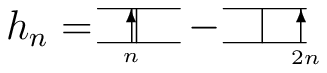}
\end{center}
\end{figure}
\noindent The fundamental weights are
$w_{i}=a_1+a_2+\cdots+a_i$ since

\begin{figure}[!h]
\begin{center}
\includegraphics[width=5in]{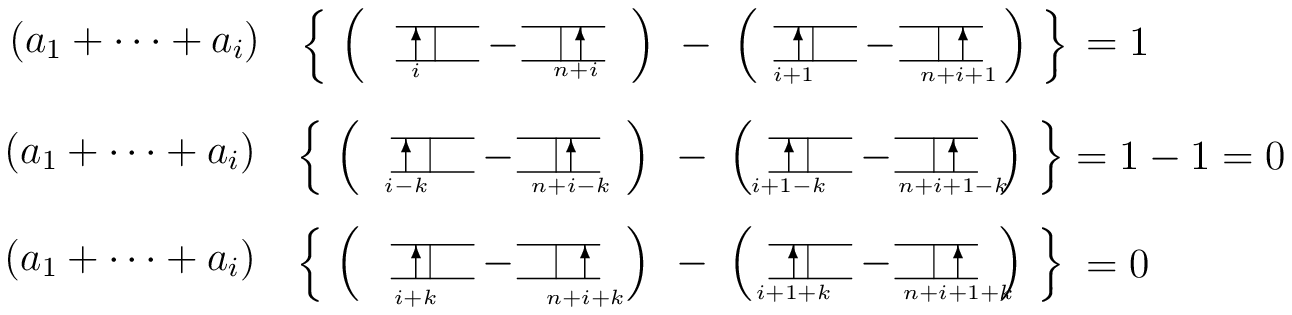}
\end{center}
\end{figure}

\subsection{Killing form of ${\mathfrak {sp}}_{2n}(\mathbb{C})$}

Let $x$ and $y$ be in ${\mathfrak h}$, we have

\begin{eqnarray*}
\langle x,y\rangle & = & \sum_{\alpha\in \Phi}\alpha(x)\alpha(y)\\
\mbox{ }& = &\sum_{i\neq j}(x_i-x_j)(y_i-y_j)+\sum_{i\neq
j}(x_i+x_j)(y_i+y_j)+
2\sum_{i}(2x_i)(2y_i)\\
\mbox{}&=&4(n+1)\sum x_i y_i\\
\mbox{}&=&4(n+1)\tr(xy).
\end{eqnarray*}

\subsection{Weyl group of ${\mathfrak {sp}}_{2n}(\mathbb{C})$}

Consider the fundamental roots of the form $\alpha_i= a_i-a_{i+1}$.
Similarly to the ${\mathfrak {sl}}_{n}(\mathbb{C})$, it is easy to
check that they generate a copy of $S_n$. Let us compute the
reflection associated to the root $\alpha_n=2a_n$. Given $h\in
{\mathfrak h}$, we have that
$S_{\alpha_{n}}(h)=h-\alpha_{n}(h)h_{\alpha_n}$, where
$h_{\alpha_n}$ is the coroot associated to the root $\alpha_n$

\begin{figure}[ht]
\begin{center}
\includegraphics[width=3.7in]{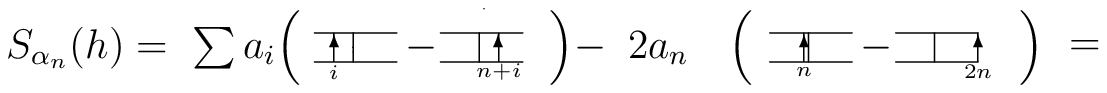}
\includegraphics[width=3.7in]{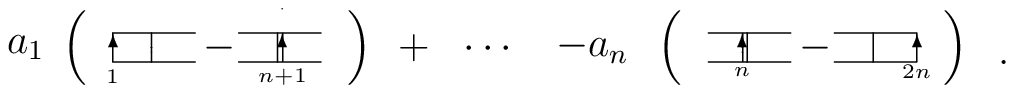}
\end{center}
\end{figure}
\noindent This reflections are the sign changes and they generate a
copy of the group $\mathbb{Z}_{2}^{n}$. Altogether the Weyl group
associated to  ${\mathfrak {sp}}_{2n}(\mathbb{C})$ is
$$C_{n}=\mathbb{Z}_{2}^{n}\rtimes S_n.$$

\subsection{ Cartan matrix and Dynkin diagram of ${\mathfrak {sp}}_{2n}(\mathbb{C})$}

$$C_n=\left( \begin{array}{cccccc}
  2 &\!\!\!\! -1 & 0 & 0 & \ldots & 0 \\
  \!\!\!\!-1 & 2 & \!\!\!\!-1 & 0 & \ldots  & 0 \\
  \vdots & \vdots & \ddots  & & \vdots & \vdots  \\
  \vdots & \vdots &  & \ddots &  & \vdots  \\
  0 & 0  & \ldots & \!\!\!\!-1 & 2 & \!\!\!\!-1 \\
  0& 0  & \ldots & 0 & \!\!\!\!-2 & 2
\end{array}\right)$$

\noindent There are $n$ vertices in this case, one for each
fundamental root. The Killing form is $\langle \alpha_i,
\alpha_{i+1}\rangle=1$, if $i=1, \dots , n-1$  and $\langle
\alpha_{n-1}, \alpha_{n}\rangle=2$. Moreover $\langle \alpha_{n-1},
\alpha_{n-1}\rangle<\langle \alpha_{n}, \alpha_{n} \rangle$, and
thus the Dynkin diagram of ${\mathfrak {sp}}_{2n}(\mathbb{C})$ has
the form
\begin{figure}[ht]
\begin{center}
\includegraphics[width=2in]{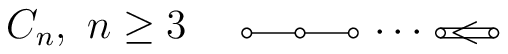}
\end{center}
\end{figure}

\subsection{Invariant functions under the action of $C_{n}=\mathbb{Z}_{2}^{n}\rtimes
S_{n}$} Let us recall that the group structure on
$\mathbb{Z}_{2}^{n}\rtimes S_n$ is given by
$$(a_,\pi)(b,\sigma)=(a\cdot\pi(b),\pi\circ\sigma)$$
where $(\pi b)_i=b_{\pi^{-1}(i)}$.

\begin{prop}\label{ACC}
$\mathbb{Z}_{2}^{n}\rtimes S_n$ acts on  $\mathbb{R}^{n}$ as
follows
$$\begin{array}{ccc}
  \mathbb{Z}_{2}^{n}\rtimes S_n \times \mathbb{R}^{n} & \longrightarrow & \mathbb{R}^{n} \\
  ((a,\pi)x) & \longmapsto & ((a,\pi)x)_i=a_ix_{\pi^{-1}(i)}
\end{array}$$
\end{prop}

\noindent Consider the polynomials
$$f_i=x_1^{2i}+x_2^{2i}+\cdots+x_n^{2i}, \ \  1\leq i\leq n$$
Each polynomial $f_i$ is invariant under the action of
$(\mathbb{Z}_{2})^{n}\rtimes S_n$ given by
$$(f(a,\pi))(x)=f((a,\pi)x).$$
The set of invariants
$$\begin{array}{ccc}
  f_1 & = & x_1^2 + x_2^2+ \cdots + x_n^2\\
  f_2 & = &  x_1^4 + x_2^4+ \cdots + x_n^4 \\
  \vdots & \mbox{ } & \mbox{ } \\
  f_n & = &  x_1^{2n}+x_2^{2n}+\cdots +x_n^{2n}
\end{array}$$
is a  basic set. This follows from the Jacobian criterion since
$$\gr(f_1)\gr(f_2)\cdots \gr(f_n)=2\cdot4\cdot6\cdots 2n= 2^n n!=|\mathbb{Z}_{2}^{n}\rtimes S_n|$$
and
$$J=2^n n!\   x_1\cdots x_n \prod_{1\leq i<j\leq n}(x_j^2-x_i^2)\neq 0.$$

\section{Orthogonal Lie Algebra ${\mathfrak
{so}}_{2n}(\mathbb{C})$}\label{sop}

Recall that the $2n$-orthogonal Lie Algebra is defined as follows
$${\mathfrak {so}}_{2n}(\mathbb{C})=\{X:X^{t}S+SX=0\}$$ where $S\in M_{2n}(\mathbb{C})$
is the matrix
$$S=\left(\begin{array}{cc}
  0 & I_n \\
  I_n & 0
\end{array} \right)$$
Eexplicitly
$${\mathfrak {so}}_{2n}(\mathbb{C})=\Bigg \{\left(
\begin{array}{cc}
  A & B \\
  C & -A^{t}
\end{array}\right);\  A,B,C\in M_{n}(\mathbb{C})\ \mbox{y}\  B=-B^{t},\  C=-C^{t} \Bigg \}.$$
${\mathfrak {so}}_{2n}(\mathbb{C})$ as a subspace of
$\textbf{Digraph}^{1}(2n,2n)$ is given by

\begin{figure}[!h]
\begin{center}
\includegraphics[width=6in]{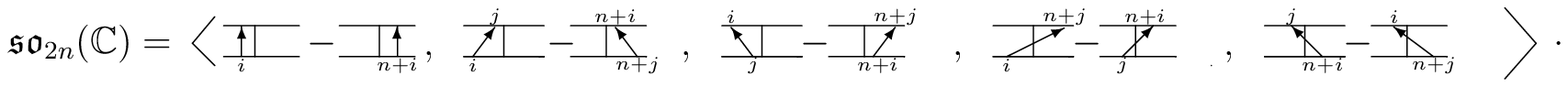}
\end{center}
\end{figure}

\noindent where $1\leq i<j\leq n$. We fix ${\mathfrak h}$ the Cartan
subalgebra of ${\mathfrak {so}}_{2n}(\mathbb{C})$ to be

\begin{figure}[!h]
\begin{center}
\includegraphics[width=3in]{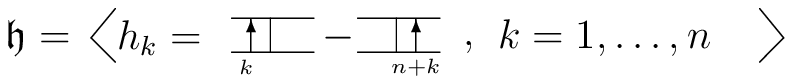}
\end{center}
\end{figure}

\newpage

\subsection{Root System of ${\mathfrak {so}}_{2n}(\mathbb{C})$}

Let $h\in {\mathfrak h}$

\begin{figure}[!h]
\begin{center}
\includegraphics[width=1.8in]{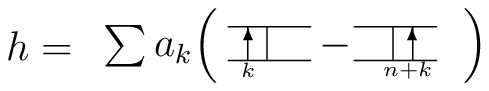}
\end{center}
\end{figure}

\noindent where $\{a_i\}$ is the base of ${\mathfrak h}^{\ast}$ dual
to the natural base of ${\mathfrak h}$. As for the case of the
symplectic algebra we define a map $T: \textbf{Digraph}^{1}(2n,2n)
\longrightarrow \textbf{Digraph}^{1}(2n,2n)$. $T$ sends a given
graph to its opposite if it does not cross the vertical line, and to
minus its opposite if it crosses the vertical line. We have again
that $T(ab)=T(b)T(a)$ For example,

\begin{figure}[!h]
\begin{center}
\includegraphics[width=1.5in]{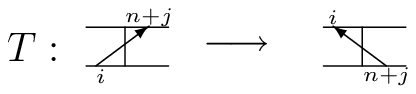}
\end{center}
\end{figure}

\noindent Let us find out the positive roots

\begin{figure}[!h]
\begin{center}
\includegraphics[width=3in]{rson.eps}
\includegraphics[width=1.7in]{rson1.eps}
\includegraphics[width=3in]{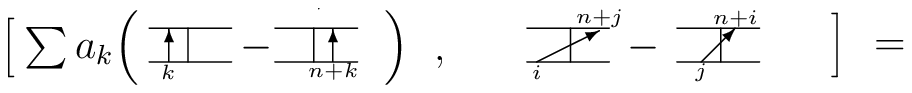}
\includegraphics[width=1.7in]{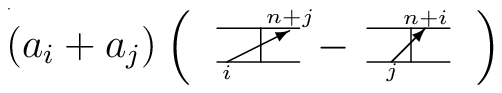}
\end{center}
\end{figure}

\noindent To get the negative roots it is enough to apply $T$ to the
positive roots. Therefore the root system is $\Phi=\{ a_i-a_j,
a_j-a_i, a_i+a_j, -a_i-a_j, \ \ 1\leq i<j\leq n\}$ and the
fundamental roots can be taken to be $\Pi=\{\alpha_1,
\dots,\alpha_n\}$ where $\alpha_i=a_i-a_{i+1}, \ \ i=1,\dots,n-1$\ \
and\ \ $\alpha_n=a_{n-1}+a_{n}$

\subsection{Coroots and weights of ${\mathfrak {so}}_{2n}(\mathbb{C})$}

\begin{enumerate}
\item{Coroot associated to the root $a_i-a_j$
\begin{figure}[!h]
\begin{center}
\includegraphics[width=2.5in]{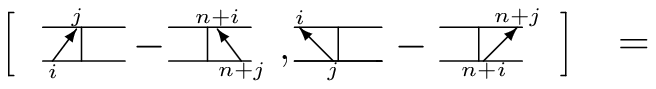}
\includegraphics[width=3.2in]{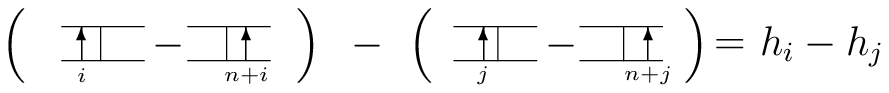}
\end{center}
\end{figure}

in this case $(a_i-a_j)(h_i-h_j)=2$, thus $h_i-h_j$ is the coroot
associated to the root $a_i-a_j$.}

\item{Coroot associated to the root $a_i+a_j$
\newpage
\begin{figure}[ht]
\begin{center}
\includegraphics[width=2.5in]{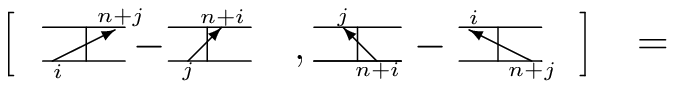}
\includegraphics[width=3.2in]{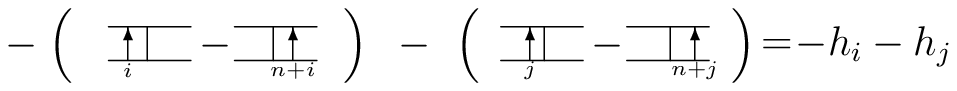}
\end{center}
\end{figure}
since $(a_i+a_j)(-h_i-h_j)=-2$, $h_i+h_j$ is the coroot associated
to the root $a_i+a_j$.}
\end{enumerate}
We concluded that  $\Phi_{c}=\{h_i-h_j, h_j-h_i, h_i+h_j,
-h_i-h_j, 1\leq i<j\leq n\}$ is the coroot system of ${\mathfrak
{so}}_{2n}(\mathbb{C})$. The set of fundamental coroots is given
by $\Pi_{c}=\{h_i-h_{i+1}, h_{n-1}+h_{n};\ \ 1\leq i\leq n-1\}$
where

\begin{figure}[!h]
\begin{center}
\includegraphics[width=3.2in]{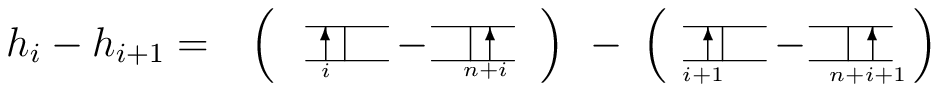}
\end{center}
\end{figure}

\noindent and

\begin{figure}[!h]
\begin{center}
\includegraphics[width=3.4in]{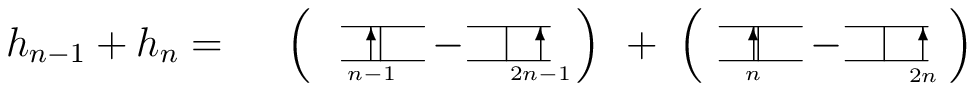}
\end{center}
\end{figure}

\noindent The fundamental weights are given by $w_i=a_1+\cdots+a_i,
\ \ i=1,\cdots,n-1$ and $
w_{n}=\frac{a_1+a_2+\cdots+a_n}{2}$ In a similar as for $\mathfrak
{sp}_{2n}(\mathbb{C})$ one can prove that
$w_{i}(h_{j}-h_{j+1})=\delta_{ij}$. For $w_{n}$ we get

\begin{figure}[!h]
\begin{center}
\includegraphics[width=4in]{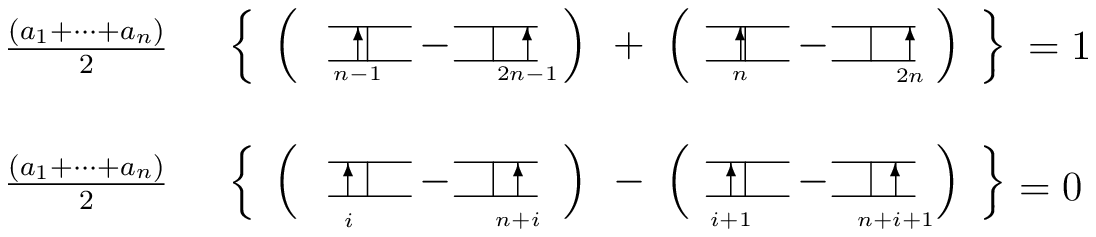}
\end{center}
\end{figure}

\subsection{Killing form of ${\mathfrak {so}}_{2n}(\mathbb{C})$}
Let $x$ and $y$ be in ${\mathfrak h}$

\begin{eqnarray*}
\langle x,y\rangle & = & \sum_{\alpha\in \Phi}\alpha(x)\alpha(y)\\
\mbox{ }& = &\sum_{i\leq j}(x_i-x_j)(y_i-y_j)+\sum_{i\leq
j}(x_i+x_j)(y_i+y_j)+
\sum_{j\leq i}(x_i-x_j)(y_i-y_j)\\
\mbox{} &\mbox{}&+\sum_{j\leq i}(x_i+x_j)(y_i+y_j)\\
\mbox{}&=&\sum_{i\neq j}2x_iy_j+2x_iy_j\\
\mbox{}&=&4(n-1)\sum x_i y_i\\
\mbox{}&=&4(n-1)\tr(xy).
\end{eqnarray*}

\subsection{Weyl group of ${\mathfrak {so}}_{2n}(\mathbb{C})$}

Consider the fundamental roots $\alpha_i= a_i-a_{i+1}$. Just as for
${\mathfrak {sl}}_{n}(\mathbb{C})$, the associated reflections
associated to these roots generate the group $S_n$. We compute the
reflections associated to the roots $\alpha_n=a_{n+1}+a_n$. Given
$h\in {\mathfrak h}$, we have
$S_{\alpha_{n}}(h)=h-\alpha_{n}(h)h_{\alpha_n}$ where $h_{\alpha_n}$
is the coroot associated to the root $\alpha_n$

\begin{figure}[!h]
\begin{center}
\includegraphics[width=5.5in]{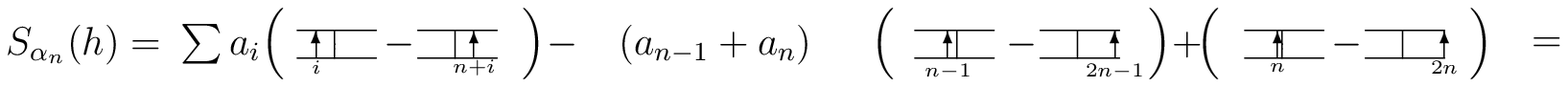}
\includegraphics[width=5in]{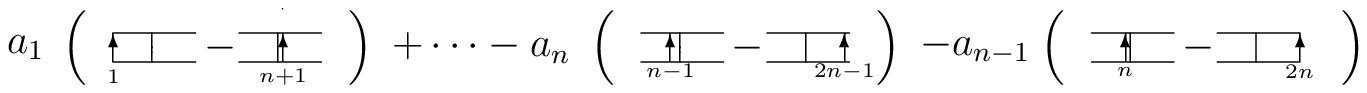}
\end{center}
\end{figure}

\noindent This reflection correspond to a change of sign. Thus we
have that the Weyl group associated with ${\mathfrak
{so}}_{2n}(\mathbb{C})$ is
$$D_{n}=\mathbb{Z}_{2}^{n-1}\rtimes S_n$$

\subsection{Cartan matrix and Dynkin diagram of ${\mathfrak {so}}_{2n}(\mathbb{C})$.}

$$D_n=\left( \begin{array}{cccccc}
  2 &\!\!\!\! -1 & 0 & 0 & \ldots&  0 \\
  \!\!\!\!-1 & 2 & \!\!\!\!-1 & 0 & \ldots  & 0 \\
  \vdots & \vdots &\ddots   & & & \vdots  \\
  0 & 0  & \ldots & 2 & \!\!\!\!-1 & \!\!\!\!-1 \\
  0 & 0  & \ldots & \!\!\!\!-1 & 2 & 0 \\
  0& 0  & \ldots & \!\!\!\!-1 & 0 & 2
\end{array}\right)$$

\noindent The Dynkin diagram has $n$ vertices  corresponding with
the fundamental roots. The Killing has the form $\langle \alpha_i,
\alpha_{i+1}\rangle=1$, if $i=1, \dots , n-2$, $\langle \alpha_{n-2}, \alpha_{n-1}\rangle=1$,
$\langle\alpha_{n-2}, \alpha_{n}\rangle=1$ and $\langle
\alpha_{n-1}, \alpha_{n}\rangle=0$.
Thus the Dynkin diagram ${\mathfrak {so}}_{2n}(\mathbb{C})$

\begin{figure}[ht]
\begin{center}
\includegraphics[width=2in]{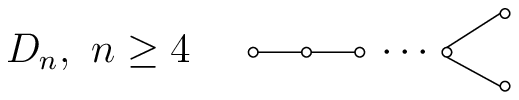}
\end{center}
\end{figure}

\subsection{Invariant functions under the action of $\mathbb{Z}_{2}^{n-1}\rtimes S_n$}

Consider the polynomials
$$\begin{array}{cccc}
  f_i & = & \displaystyle {\sum_{j=1}^{n}x_j^{2i}}, & 1\leq i \leq n-1 \\
  f_n & = & x_1\cdots x_n & \mbox{ }
\end{array}$$
clearly each $f_i$ is invariant under the action of
$\mathbb{Z}_{2}^{n-1}\rtimes S_n$. It is easy to check that
$$\gr(f_1)\gr(f_2)\dots \gr(f_n)= 2^{n-1} n!=|\mathbb{Z}_{2}^{n-1}\rtimes S_n|$$
and
$$J=(-2)^{n-1}(n-1)! \displaystyle {\prod_{1\leq i<j\leq n}(x_j^2-x_i^2)} \neq 0,$$
so the Jacobian criterion  tell us that ${f_1,\dots, f_n}$ is a
basic set of invariants.

\section{Orthogonal algebra ${\mathfrak
{so}}_{2n+1}(\mathbb{C})$}\label{soi}

We orthogonal odd algebra
$${\mathfrak {so}}_{2n+1}(\mathbb{C})=\{X:X^{t}S+SX=0\}$$ where $S\in
M_{2n+1}(\mathbb{C})$ is of the form
$$S=\left(\begin{array}{ccc}
  0 & I_n & 0\\
  I_n & 0 & 0\\
  0 & 0 & 1
\end{array} \right)$$
In an explicit form
$${\mathfrak {so}}_{2n+1}(\mathbb{C})=\Bigg\{ \left(
\begin{array}{ccc}
  A & B &-H^{t}\\
  C & -A^{t}& -G^{t}\\
  G & H & 0
\end{array}\right);\   A,B,C\in M_{n}(\mathbb{C}), H,G\in M_{1\times n}(\mathbb{C}), 0\in\mathbb{C}
\ \ \mbox{y} \  B=-B^{t},\  C=-C^{t} \Bigg\}$$ ${\mathfrak
{so}}_{2n+1}(\mathbb{C})$ as a subspace of
$\textbf{Digraph}^{1}(2n+1,2n+1)$ is given by

\begin{figure}[!h]
\begin{center}
\includegraphics[width=6in]{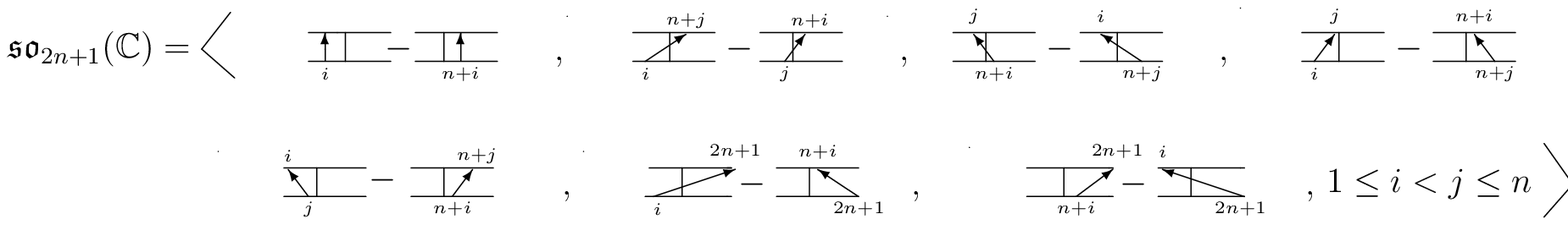}
\end{center}
\end{figure}

\noindent Let us fix ${\mathfrak h}$ a Cartan subalgebra of
${\mathfrak {so}}_{2n+1}(\mathbb{C})$
\begin{figure}[ht]
\begin{center}
\includegraphics[width=2.5in]{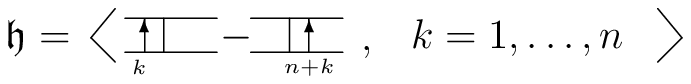}
\end{center}
\end{figure}

\subsection{Root system of ${\mathfrak {so}}_{2n+1}(\mathbb{C})$}

Let $h\in {\mathfrak h}$,

\vspace{-0.5cm}

\begin{figure}[ht]
\begin{center}
\includegraphics[width=1.7in]{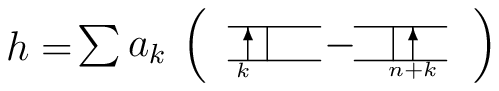}
\end{center}
\end{figure}

\noindent where $\{a_i\}$ is a base of ${\mathfrak h}^{\ast}$ dual to the
natural base of ${\mathfrak h}$. We compute the positive roots. The
negative roots are obtain applying the following antimorphism to the
positive roots.\\
$T: \textbf{Digraph}^{1}(2n+1,2n+1)
\longrightarrow \textbf{Digraph}^{1}(2n+1,2n+1)$. For example

\begin{figure}[!h]
\begin{center}
\includegraphics[width=1.5in]{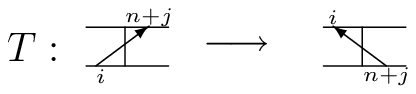}
\end{center}
\end{figure}

\begin{figure}[!h]
\begin{center}
\includegraphics[width=3.2in]{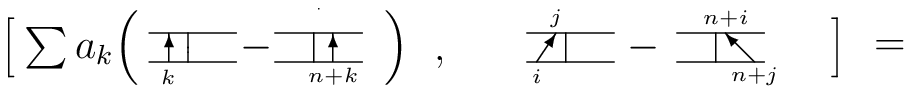}
\includegraphics[width=1.8in]{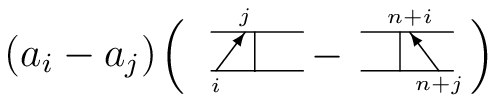}
\includegraphics[width=3.2in]{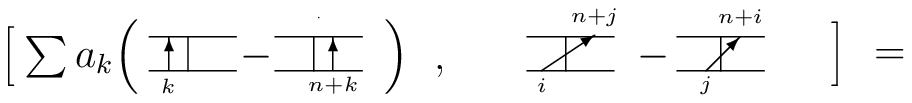}
\includegraphics[width=1.8in]{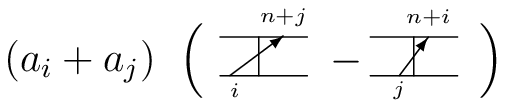}
\includegraphics[width=3.2in]{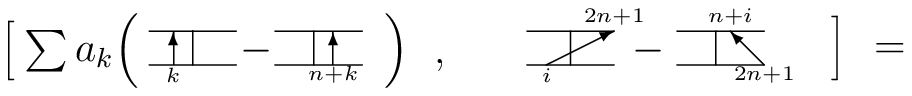}
\includegraphics[width=1.5in]{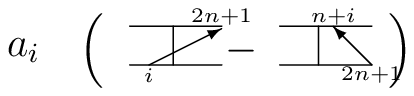}
\end{center}
\end{figure}

\noindent Thus the set of roots is $\Phi=\{a_i-a_j,
a_j-a_i,a_i+a_j,a_i,-a_i\}.$ The fundamental roots are
$\Pi=\{\alpha_1,\dots,\alpha_n\}$ where $\alpha_i=a_i-a_{i+1}, \ \
i=1,\dots,n-1$\ \ and\ \ $\alpha_n=a_{n}$. In pictures the root
system of ${\mathfrak {so}}_{7}(\mathbb{C})$ looks like
\begin{figure}[!h]
\begin{center}
\includegraphics[height=2cm]{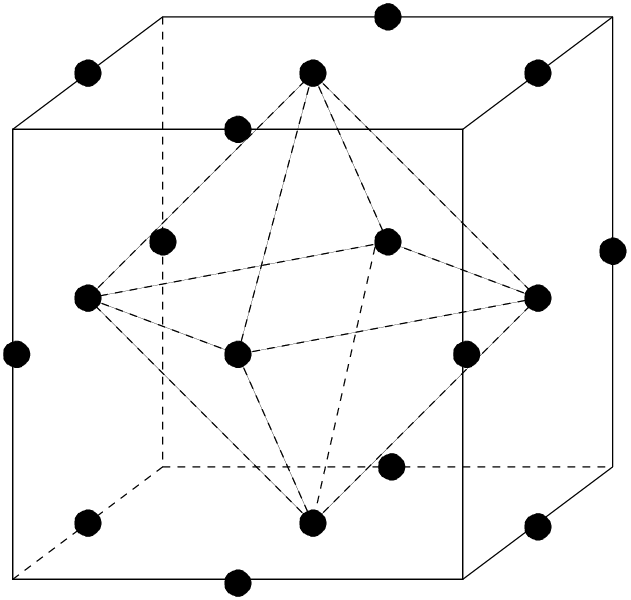}
\end{center}
\end{figure}

\subsection{Coroots and weights of ${\mathfrak {so}}_{2n+1}(\mathbb{C})$}

\begin{enumerate}
\item{Coroots associated to the root $a_i-a_j$

\begin{figure}[!h]
\begin{center}
\includegraphics[width=2.5in]{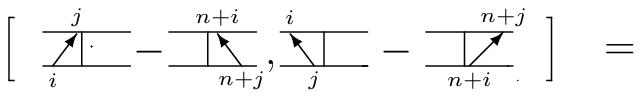}
\includegraphics[width=3in]{coson1.eps}
\end{center}
\end{figure}

Since $(a_i-a_j)(h_i-h_j)=2$, we see that  $h_i-h_j$ is the coroot
associated to the root $a_i-a_j$.}

\item{Coroot associated to the root $a_i+a_j$

\begin{figure}[!h]
\begin{center}
\includegraphics[width=2.5in]{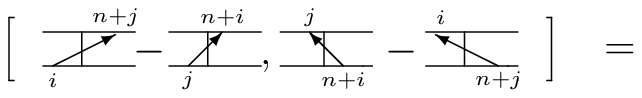}
\includegraphics[width=3.2in]{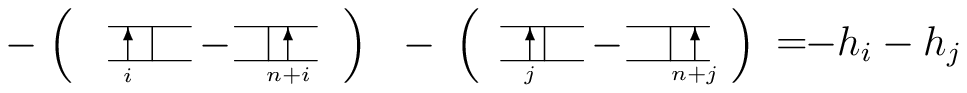}
\end{center}
\end{figure}

Here $(a_i+a_j)(-h_i-h_j)=-2$, and thus $h_i+h_j$ is the coroot
associated to the root $a_i+a_j$.}

\item{Coroot associated to the root $a_i$

\begin{figure}[!h]
\begin{center}
\includegraphics[width=2.5in]{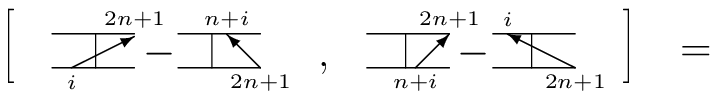}
\includegraphics[width=2 in]{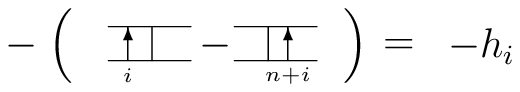}
\end{center}
\end{figure}

Thus $2h_i$ is the coroot associated to the root $a_i$.}
\end{enumerate}
We have that  $\Phi_{c}=\{h_i-h_j, h_j-h_i, h_i+h_j, -h_i-h_j,
2h_i,-2h_i, \ 1\leq i<j\leq n \}$ is the coroot system of
${\mathfrak {so}}_{2n+1}(\mathbb{C})$. The set of fundamental coroot
has the form $\Pi_{c}=\{h_{i}-h_{i-1}, 2h_{n};\ \ 1\leq i\leq n-1\}$
where

\begin{figure}[!h]
\begin{center}
\includegraphics[width=3.4in]{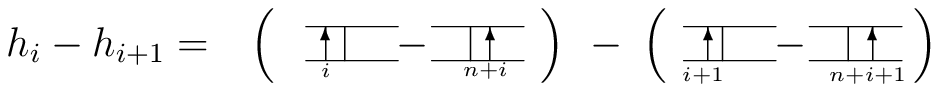}
\end{center}
\end{figure}
\noindent and

\begin{figure}[!h]
\begin{center}
\includegraphics[width=3.4in]{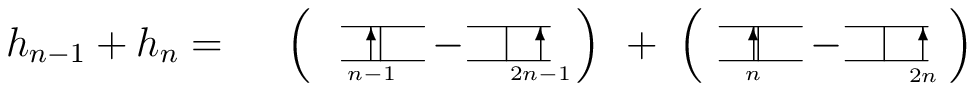}
\end{center}
\end{figure}

\noindent The fundamental weights are $w_{i}=a_1+\cdots+a_{i},\ \
i=1,\cdots,n-1$ and ${\displaystyle
w_{n}=\frac{a_1+\cdots+a_n}{2}.}$ In a similar fashion to the
${\mathfrak {so}}_{2n}(\mathbb{C})$ case we have that
$w_{i}(h_j-h_{j+1})=\delta_{ij}$ and $w_{i}(2h_n)=\delta_{in}$.

\subsection{Killing form of ${\mathfrak {so}}_{2n+1}(\mathbb{C})$}
Given $x$ and $y$ in ${\mathfrak h}$

\begin{eqnarray*}
\langle x,y\rangle & = & \sum_{\alpha\in \Phi}\alpha(x)\alpha(y)\\
\mbox{ }& = &\sum_{i\leq j}(x_i-x_j)(y_i-y_j)+\sum_{i\leq
j}(x_i+x_j)(y_i+y_j)+
\sum_{j\leq i}(x_i-x_j)(y_i-y_j)\\
\mbox{} &\mbox{}&+\sum_{j\leq i}(x_i+x_j)(y_i+y_j)+\sum_{i}x_i y_i+\sum_{ i}(2x_i)(2y_i)\\
\mbox{}&=&(4n-2)\sum x_i y_i\\
\mbox{}&=&(4n-2)\tr(xy)
\end{eqnarray*}

\subsection{Weyl group of ${\mathfrak {so}}_{2n+1}(\mathbb{C})$}

Consider the fundamental roots $\alpha_i= a_i-a_{i+1}$. Just like
for ${\mathfrak {sl}}_{n+1}(\mathbb{C})$, the reflections associated
to these roots generate the symmetric group $S_n$. Let us analyze
the reflection associated to the root $\alpha_n=a_n$. Let $h\in
{\mathfrak h}$, we have
$S_{\alpha_{n}}(h)=h-\alpha_{n}(h)h_{\alpha_n}$ where $h_{\alpha_n}$
is the coroot associated to the root $\alpha_n$

\begin{figure}[!h]
\begin{center}
\includegraphics[width=4in]{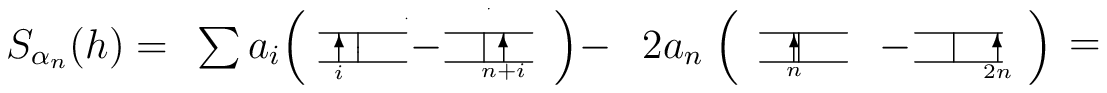}
\includegraphics[width=3.3in]{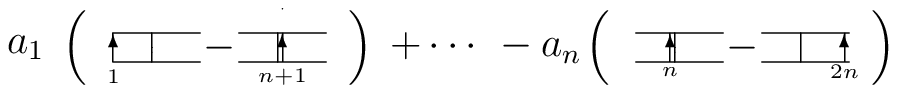}
\end{center}
\end{figure}

\noindent This reflections represent sign changes and generate the
group $\mathbb{Z}_{2}^{n}$, therefore the Weyl group associated with
${\mathfrak {so}}_{2n+1}(\mathbb{C})$ is
$$B_n=\mathbb{Z}_{2}^{n}\rtimes S_n$$

\subsection{Cartan matrix and Dynkin diagram of ${\mathfrak {so}}_{2n+1}(\mathbb{C})$.}

$$B_n=\left( \begin{array}{cccccc}
  2 &\!\!\!\! -1 & 0 & 0 & \ldots &  0 \\
  \!\!\!\!-1 & 2 &\!\!\!\! -1 & 0 & \ldots &  0 \\
  \vdots & \vdots & \ddots &  &  & \vdots  \\
  \vdots & \vdots &  &\ddots  &  & \vdots  \\
  0 & 0 & \ldots &\!\!\!\! -1 & 2 & \!\!\!\!-2 \\
   0 & 0 & \ldots & 0 & \!\!\!\!-1 & 2
\end{array} \right)$$
The diagram has $n$ vertices, one for each fundamental root. The
killing form is given by $\langle \alpha_i, \alpha_{i+1}\rangle=1$,
if $i=1, \dots , n-1$ and $\langle \alpha_{n-1},
\alpha_{n}\rangle=2$. Furthermore $\langle \alpha_{n-1},
\alpha_{n-1}\rangle>\langle \alpha_{n}, \alpha_{n} \rangle$, and
thus, the Dynkin diagram of ${\mathfrak {so}}_{2n+1}(\mathbb{C})$
has form
\begin{figure}[ht]
\begin{center}
\includegraphics[width=2in]{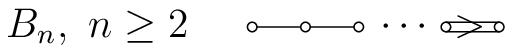}
\end{center}
\end{figure}

\subsection*{Acknowledgment}
Thanks to Manuel Maia for helping us with Latex.

$$\begin{array}{c}
\mbox{Rafael D\'\i az. Instituto Venezolano de Investigaciones Cient\'\i ficas (IVIC).} \ \  \mbox{\texttt{radiaz@ivic.ve}} \\
\!\!\!\!\!\mbox{Eddy Pariguan. Universidad central de Venezuela (UCV).} \ \  \mbox{\texttt{eddyp@euler.ciens.ucv.ve}} \\
\end{array}$$

\end{document}